\RequirePackage{fix-cm}
\documentclass[smallextended, envcountsame]{svjour3}
\smartqed
\def\makeheadbox{\relax}
\journalname{}

\usepackage[unicode = false,
    pdftoolbar = true,
    colorlinks = true,
    linkcolor = blue,
    citecolor = blue,
    filecolor = black,
    urlcolor = blue,
    breaklinks = true]{hyperref}
\usepackage[ruled]{algorithm2e}
\usepackage{float}
\usepackage{geometry}
\usepackage[utf8]{inputenc}
\usepackage{listings}           
\usepackage{amsmath,amssymb}            
\usepackage{xcolor}             
\usepackage{graphicx}           
\usepackage[ruled]{algorithm2e}
\usepackage{verbatim}
\usepackage{geometry}
\usepackage{enumerate}
\usepackage[thinlines, thiklines]{easybmat}
\usepackage{hhline}
\usepackage{nicematrix}
\usepackage{mathtools,nccmath}
\usepackage[noblocks]{authblk}
\usepackage{caption}
\usepackage[section]{placeins}

\DeclareMathOperator*{\Id}{Id}

\DeclareMathOperator{\diag}{Diag}

\DeclareMathOperator{\Diag}{Diag}

\newcommand{\R}{\mathbb{R}}

\newcommand{\bbm}{\begin{bmatrix}}
\newcommand{\ebm}{\end{bmatrix}}

\newcommand{\zero}{\mathbf{0}}

\newcommand{\Span}{\operatorname{span}}
\newcommand{\sh}{\nabla^2_s}

\newcommand{\Proj}{\operatorname{P}}
\newcommand{\dom}{\operatorname{dom}}

\newcommand{\N}{\operatorname{\mathbb{N}}}

\newcommand{\U}{\operatorname{\overline{U}}}

\newcommand{\dnsf}{\delta_s^2}
\newcommand{\epsfS}{\varepsilon f(x^0;S)}

\newcommand{\ns}{\nabla_s}

\newcommand{\cshtwo}{\nabla^2_{c}}

\newcommand\scalemath[2]{\scalebox{#1}{\mbox{\ensuremath{\displaystyle #2}}}}

\newcommand{\shessti}{\nabla_s^2 f(x^0;S,T_{1:m})}
\newcommand{\cshesstitwo}{\nabla_{c}^2 f(x^0;S,T_{1:m})}

\newcommand{\Projuv}{\Proj_{S, \T}}

\newcommand{\Projuivi}{\Proj_{S,T_{1:m}}}

\newcommand{\cshd}{d\nabla^2 f(x^0;S)}

\newcommand{\Ti}{T_{1:m}}
\newcommand{\T}{\overline{T}}
\newcommand{\stress}{\nabla_s^3 f(x^0;S,\T,\U)}
\newcommand{\dst}{\mathbf{\delta}_{\mathbf{s}}^\mathbf{3}}
\newcommand{\shesst}{\nabla^2_sf(x^0;S,\T)}
\newcommand{\shessv}{\nabla^2_sf(x^0;S,hv)}
\newcommand{\cshesst}{\nabla^2_cf(x^0;S,\T)}
\newcommand{\cshessv}{\nabla^2_cf(x^0;S,hv)}

\newgeometry{vmargin={35mm}, hmargin={19mm,22mm}}

\begin{document}

\title{\centering Using generalized simplex methods to approximate derivatives}
\author{\centering Gabriel Jarry-Bolduc \and Chayne Planiden}
\institute{
          Gabriel Jarry-Bolduc\at
            Department of Mathematics and Computing, Mount Royal University, \\
            Calgary, Alberta, Canada. \\
              ORCID 0000-0002-1827-8508\\
              \email{gabjarry@alumni.ubc.ca} 
              \and
              Chayne Planiden\at School of Mathematics and Applied Statistics, University of Wollongong, \\ Wollongong, NSW, 2500, Australia. Research supported by University of Wollongong.  \\ ORCID 0000-0002-0412-8445 \\
              \email chayne@uow.edu.au
              }
\date{\today}

\titlerunning{Using generalized simplex methods to approximate derivatives}

\maketitle
\begin{abstract}
This paper presents two methods for approximating a proper subset of the entries of a Hessian using only function evaluations. It is also shown how to approximate a Hessian-vector product with a minimal number of function evaluations.  These approximations are obtained  using the techniques called \emph{generalized simplex Hessian} and \emph{generalized centered simplex Hessian}.  We show how to choose the matrices of directions involved in the computation of these two techniques  depending on the entries of the Hessian of interest. We discuss the number of function evaluations required in each case and develop a general formula to approximate all order-$P$ partial derivatives. Since only function evaluations are required to compute the methods discussed in this paper, they are suitable for use in derivative-free  optimization methods.
\end{abstract}

\section{Introduction}\label{sec:intro}

Approximating Hessians  is an important topic in numerical analysis and optimization. The Hessian of a function  captures the curvature of the function, thus providing additional information that the gradient does not have and aiding in the optimization process. There exist many approaches to obtain exact Hessians or approximate  Hessians, including automatic differentiation \cite{griewank1991}, graph coloring approach \cite{coleman1984estimation,gebremedhin2005color}, Lagrange polynomials \cite{conn2009introduction},  Newton fundamental polynomials \cite{conn2009introduction}, regression nonlinear models or underdetermined interpolating models \cite{conn2008bgeometry}. Arguably, two of the most well-known methods to approximate Hessians are \emph{forward-finite-difference approximation} and \emph{centered-finite-difference} approximation \cite[Section 4.6]{andrei2021derivative}. 

In derivative-free optimization (DFO) methods, it is usually assumed that derivative information is not directly available to the individual conducting the optimization process. For this reason,  true Hessians are not employed.  Approximate Hessians have been used in DFO methods since at least 1970 \cite{Winfield1970}.  Researchers from the DFO community have previously explored methods to approximate full Hessians or some of the entries of the Hessian.  In \cite{Custodio2007}, the authors outline an idea for a \emph{simplex Hessian} that is constructed via quadratic interpolation through $(n+1)(n+2)/2$ well-poised sample
points.  They further posit that if only the diagonal entries are desired, then $2n+1$ sample points are sufficient. These ideas are formalized in \cite{conn2008bgeometry} through quadratic interpolation and analyzed through the use of Lagrange polynomials. Obtaining an approximation of the diagonal component of a Hessian is also discussed in \cite{coope2021gradient,jarry2022approximating}, and in \cite{hare2020error} it is shown that the diagonal entries can be obtained for free (in terms of function evaluations) if the gradient has been previously approximated  via the \emph{(generalized) centered simplex gradient} technique.

Several  DFO algorithms employ approximate Hessians to solve  optimization problems \cite{custodio2010incorporating,Kelley2011,powell1998use,powell2004least,powell2004use,powell2004updating,powell2006newuoa,powell2007view,powell2008,wild2008}. 
To  develop strong  convergence results, DFO algorithms   rely on  techniques to  approximate Hessians and gradients in a manner that has controllable error bounds. In other words, the accuracy of the approximate Hessians and gradients can be controlled. In \cite{hare2023hessianpublished,jarry2023numerical},  two techniques based on simple matrix algebra to approximate a full Hessian, called the \emph{generalized simplex Hessian (GSH)} and the \emph{generalized centered simplex Hessian (GCSH)}, are introduced.  It is shown that the GSH is an order-1 accurate approximation of the full Hessian and that the GCSH is an order-2 accurate approximation of the full  Hessian.  The GSH  can be viewed as a generalization of the \emph{simplex Hessian} discussed in  \cite{conn2009introduction,Custodio2007}.  The simplex Hessian requires  $(n+1)(n+2)/2$  sample points poised for quadratic interpolation. On the other hand, the GSH and the GCSH are well-defined as long as the matrices of directions utilized are nonempty. Hence, they offer enough flexibility to approximate a proper subset of the entries of a Hessian. 

  In many optimization problems, approximating  the Full Hessian  may be too expensive. For this reason, only approximating the appropriate entries of the Hessian may be necessary  to develop efficient methods. In  \cite{mohammad2018structured}, an approximation of the diagonal entries of the Hessian is used to solve large-scale nonlinear least-square problems.   In \cite{andrei2020diagonal}, it is shown that a  quasi-Newton method relying on an approximation of the diagonal entries of the Hessian is the top performer against well-established algorithms.   In \cite{akrotirianakis2004role}, the off-diagonal entries of the Hessian are utilized to construct tight convex underestimators for nonconvex functions. In presence of a sparse Hessian, approximating  some of the rows/columns of the Hessian may be sufficient to capture  enough information  about the curvature of the  function. 

  Another interesting topic related to the one of approximating a subset of the entries of a Hessian is to approximate a Hessian-vector product.  A technique that does not require to store the full Hessian was introduced by Pearlmutter in \cite{pearlmutter1994fast}.  Hessian-vector products appear  in Newton-based methods  for nonlinear optimization. These methods require the solution of a linear system  in which the matrix is the Hessian of a function and the vector is the negative gradient. Iterative methods are usually employed to find a solution \cite{song2021modeling}.  This can be done by approximating  the Hessian-vector product  with finite-differences \cite{nocedal2006numerical} or by using automatic differentiation technique as described in \cite{hicken2014inexact}.  
  
The  main goal of this paper is to investigate how to approximate a proper subset of the entries of the Hessian and a Hessian-vector product with the  GSH and the GCSH. Error bounds are provided, showing that the GSH can provide  order-1 accuracy of the appropriate subset of the entries of the Hessian. The results show that the GSH is a  versatile technique that can be used to  to do ``everything'' when approximating Hessians. Using the GCSH,  error bounds show that we can obtain order-2 accuracy on the appropriate subset of the entries of the  Hessian. A secondary goal is to introduce a general recursive formula to approximate all order-$P$ partial derivatives. Optimization methods relying on high-order (i.e., higher than two) derivatives  have recently gained traction  (for example, see \cite{cartis2020sharp,grapiglia2020tensor,nesterov2023inexact}).

This paper is organized as follows. Section \ref{sec:prel} contains a description of the notation and some needed definitions, including those of the generalized simplex gradient (GSG), the GSH and the GCSH. 
In Section \ref{sec:partialHessian}, we investigate how to approximate a proper subset of the entries of a Hessian. Diagonal entries and the relation between the GSH and the \emph{centered simplex Hessian diagonal} are discussed in Section \ref{sec:partialHessian} as well, with details provided on how to approximate the off-diagonal entries of a Hessian and a row/column of a Hessian. Error bounds are provided in each section. The properties of the matrices of directions, number of function evaluations required, and error bounds are described. In Section \ref{sec:hvp}, how to approximate a Hessian-vector product is discussed. 
 It is shown that a Hessian-vector product can be approximated with $2n+1$ function evaluations with the GSH and  $4n-1$ function evaluations with the GCSH. In Section \ref{sec:nablaP}, a formula to approximate all order-$P$ partial derivatives is introduced. Finally, Section \ref{sec:conclusion} contains concluding remarks and recommends areas of future research in this vein.

\section{Preliminaries}\label{sec:prel}

In this section, we clarify the notation and introduce background results necessary to understand the following sections. Most of this section is identical to \cite[Section 2]{hare2023hessianpublished}. 

Throughout this paper, we use the standard notation found in \cite{rockwets}. The domain of a function $f$ is denoted by $\dom f$. The transpose of matrix $A$ is denoted by $A^\top$. We work in  the finite-dimensional space $\R^n$ with inner product $x^\top y=\sum_{i=1}^nx_iy_i$ and induced norm $\|x\|=\sqrt{x^\top x}$. The identity matrix in $\R^{n \times n}$ is denoted by $\Id_n$. We use $e_n^i \in \R^n$ where $i \in \{1, \dots, n\},$ to denote the standard unit basis vectors in $\R^n$, i.e.\ the $i$\textsuperscript{th} column of $\Id_n.$ When there is no ambiguity about the dimension, we may omit the subscript and simply  write $e^i$ or $\Id.$ The zero vector in $\R^n$ is denoted by $\zero_n$ and the zero matrix in $\R^{n \times m}$ is denoted by $\zero_{n \times m}$. The entry in the $i$\textsuperscript{th} row and $j$\textsuperscript{th} column of a matrix $A \in \R^{n \times m}$ is denoted by $A_{i,j}.$ If the matrix  already involves a subscript, say $k$, then we use the notation $[A_k]_{i,j}$. The matrix $D=\diag(v)= \diag[v_1, \dots, v_n] \in \R^{n \times n},$ where $v \in \R^n,$  represents a diagonal matrix  with  diagonal entries $D_{j,j}=v_j$ for all $j \in \{1, \dots, n\}.$ 
The span of a matrix $A \in \R^{n \times m},$ denoted by $\Span A$, is the column space of $A.$ That is the space generated by all linear combinations of the columns in $A.$ 
The \emph{Minkowski sum} of two sets of vectors $A$ and $B$ is denoted by $A \oplus B$ and defined as follows:
\begin{align*}
    A \oplus B=\{a+b: a\in A, b \in B\}.
\end{align*}
If $A$ is a set that contains a single vector $x$, we write $x\oplus B$ rather than $\{x\}\oplus B$ or $A\oplus B$.
Given a matrix $A \in \R^{n \times m},$ we use the induced matrix norm 
\begin{align*}
    \Vert A \Vert=\Vert A \Vert_2=\max \{ \Vert Ax\Vert_2 \, : \, \Vert x \Vert_2=1 \}.
\end{align*}
We denote by $B_n(x^0;\Delta)$ and $\overline{B}_n(x^0;\Delta)$ the open and closed balls, respectively, centered at $x^0 \in \R^n$ with finite radius $\Delta>0$.

In order to introduce the definitions of the GSG, the GSH and the GCSH, we require the Moore--Penrose pseudoinverse of a matrix.

\begin{definition}[Moore-Penrose pseudoinverse] \cite[Chapter 17]{roman2007}\label{def:mpinv}
Let $A \in \R^{n \times m}$. The unique matrix $A^{\dagger}\in\R^{m \times n}$ that satisfies the following four equations is called the Moore-Penrose pseudoinverse of $A$:
\begin{align*}
&(i) \, AA^{\dagger}A=A,&&(ii) \, A^{\dagger}AA^{\dagger}=A^{\dagger},
&&& (iii) \, (AA^\dagger)^\top=AA^{\dagger}, &&&& (iv) \,(A^{\dagger}A)^\top=A^{\dagger}A.
\end{align*}
\end{definition}

The Moore--Penrose pseudoinverse $A^\dagger$ is not always an inverse of $A$, but the following two properties hold.
\begin{itemize}
\item If $A$ has full column rank $m$, then $A^\dagger$ is a left inverse of $A$. That is, $A^\dagger A=\Id_m$ and
\begin{equation} \label{eq:fullcolumnrank}
 A^\dagger=(A^\top A)^{-1}A^\top.
 \end{equation}
\item If $A$ has full row rank $n$, then $A^\dagger$ is a right inverse of $A.$ That is, $AA^\dagger=\Id_n$ and 
\begin{equation} \label{eq:fullrowrank}
A^\dagger=A^\top (AA^\top)^{-1}.
\end{equation}
\end{itemize}

Next, we provide the definition of the generalized simplex gradient (GSG). Recall that the simplex gradient is the gradient of the linear interpolation function defined by the points $x^0\oplus S$; the GSG is an extension of this idea that does not require $S$ to be square  or full rank.

\begin{definition}[Generalized simplex gradient]\cite[Definition 2.5]{hare2023hessianpublished} Let $f:\dom f \subseteq \R^n\to\R$ and let $x^0 \in \dom f$ be the point of interest.  Let $S\in \R^{n \times m}$ with $x^0 \oplus S \subset \dom f.$ The \emph{generalized simplex gradient of $f$ at $x^0$ over $S$} is denoted by $\nabla_s f(x^0;S)$ and defined by
\begin{align*} \nabla_s f(x^0;S)&=(S^\top)^\dagger \delta_s f(x^0;S) \in \R^n 
\end{align*} 
where 
\begin{align*}
\delta_s f(x^0;S)&=\begin{bmatrix} f(x^0+s^1)-f(x^0) \\ \vdots \\ f(x^0+s^m)-f(x^0) \end{bmatrix} \in \R^m.
\end{align*}
\end{definition}
 
Next, we clarify the key notation used in the construction of the GSH and the GCSH.  Within, we write a set of vectors in matrix form, by which we mean that each column of the matrix is a vector in the set. Let
\begin{align*}
    &S=\begin{bmatrix} s^1&s^2&\cdots&s^m\end{bmatrix} \in \R^{n \times m} ~~\mbox{and}~~T_j=\begin{bmatrix} t^1_j&t^2_j&\cdots&t^{k_j}_j \end{bmatrix} \in \R^{n \times k_j},  j\in\{1,\ldots,m\},
\end{align*}
be sets of directions contained in $\R^n$.  Define
    $$\Ti=\{T_1,\ldots,T_m\},$$

\noindent and
    $$\Delta_S=\max\limits_{j \,\in\{1,\ldots,m\}}\|s^j\|, \quad \Delta_{T_j}=\max\limits_{\ell \,\in\{1,\ldots,k_j\}}\|t^\ell_j\|, \quad \Delta_T=\max\limits_{j\in\{1,\ldots,m\}}\Delta_{T_j}.$$

 \noindent The \emph{normalized matrices} $\widehat{S}$ and $\widehat{T_j}$ are  respectively defined by
 \begin{align} \label{eq:hat}
\widehat{S}=\frac{1}{\Delta_S} S, \quad \widehat{T_j}=\frac{1}{\Delta_{T_j}} T_j, \quad j \in \{1, \dots, m\}.
 \end{align}

\noindent In this paper, it  is always assumed that the matrix $S$ and all matrices $T_{j}$ are non-empty and have non-null rank. This ensures that the matrices in \eqref{eq:hat} are well-defined.

\begin{definition}[Generalized simplex Hessian] \cite[Definition 3.1]{hare2023hessianpublished}\label{def:gsh}
Let $f:\dom f \subseteq \R^n \to \R$ and let $x^0 \in \dom f$ be the point of interest.  Let $S=\begin{bmatrix} s^1&s^2&\cdots &s^m \end{bmatrix}  \in \R^{n \times m}$ and $ T_j \in \R^{n \times k_j}$ with $x^0 \oplus T_j,x^0 \oplus S,x^0+s^j \oplus T_j$ contained in  $\dom f$ for all $j \in \{1, \dots, m\}.$ The \emph{ generalized simplex Hessian  of $f$ at $x^0$ over $S$ and $\Ti$} is denoted by $\sh f(x^0;S,\Ti)$ and defined by
\begin{equation*}\sh f(x^0;S,\Ti)=(S^\top)^\dagger\dnsf f(x^0;S,\Ti),  \end{equation*} 
where
\begin{equation*}
    \dnsf f(x^0;S;\Ti)=\left[\begin{array}{c}(\ns f(x^0+s^1;T_1)-\ns f(x^0;T_1))^\top\\(\ns f(x^0+s^2;T_2)-\ns f(x^0;T_2))^\top\\\vdots\\(\ns f(x^0+s^m;T_m)-\ns f(x^0;T_m))^\top\end{array}\right]\in \R^{m \times n}.
\end{equation*}
In the case $T_1=T_2=\cdots=T_m$, we use $\overline{T}=T_j \in \R^{n \times k}$ to simplify notation and write $\sh f(x^0;S,\T)$ to emphasize the special case.
\end{definition}

Note that the number of  columns $m$ in $S$ can be any  positive integer and each  number of columns $k_j$ in $T_j$ can be any positive integer for all $j \in \{1, \dots, m\}.$

A ``centered'' version of the GSH can also be defined.

\begin{definition}[Generalized centered simplex Hessian]\cite[Definition 3.2]{hare2023hessianpublished} \label{def:gcsh}
 Let $f:\dom f \subseteq \R^n \to \R$ and let $x^0 \in \dom f$ be the point of interest.  Let $S  \in \R^{n \times m}$ and $ T_j \in \R^{n \times k_j}$ with $x^0 \oplus  S  \oplus  T_j,  x^0 \oplus (-S) \oplus (-T_j),  x^0 \oplus (\pm S),$ and $x^0 \oplus (\pm T_j)$ contained   in  $\dom f$ for all $j \in \{1, \dots, m\}.$  The  \emph{generalized centered simplex Hessian  of $f$ at $x^0$ over $S$ and $\Ti$} is denoted by $\cshtwo f(x^0;S,T_{1:m})$ and defined by   
\begin{equation*}
\cshtwo f(x^0;S,T_{1:m})=\frac{1}{2} \left ( \shessti + \nabla_s^2 f(x^0;-S,-\Ti) \right ).
\end{equation*}
\end{definition}

It turns out that error bounds can be defined between the GSH (GCSH) and some of the entries of the true Hessian. The appropriate entries of the true Hessian are obtained via a projection operator (previously defined in \cite[Section 3]{hare2023hessianpublished}). 
The projection operator involves all matrices of directions utilized to  compute  the GSH (GCSH).
Given matrices $S \in \R^{n \times m}$ and  $T_j \in \R^{n \times k_j}$, the projection of a matrix $M \in\R^{n \times n}$ onto $S$  and $\Ti$ is denoted by $\Projuivi M$ and defined by
\begin{align*}
\Projuivi M&=\sum_{j=1}^m (S^\top)^\dagger e^j_m (e^j_m)^\top S^\top M T_j T_j^\dagger.
\end{align*}
In the case where $T_1=T_2=\dots=T_m=\T,$ the projection of $M$ onto $S$ and $\T$ is denoted by $\Projuv M,$ and reduces to 
\begin{align*}
\Projuv M&=\sum_{j=1}^m (S^\top)^\dagger e^j_m (e^j_m)^\top S^\top M \T \T^\dagger \\
&=(S^\top)^\dagger \left ( \sum_{j=1}^m  e^j_m (e^j_m)^\top \right )S^\top M\T\T^\dagger \\
&=(S^\top)^\dagger {\Id}_m S^\top M\T\T^\dagger \\
&=(S^\top)^\dagger S^\top M\T\T^\dagger.
\end{align*}

\noindent Note that $\Projuivi$ is a linear operator. If one of the following three conditions is satisfied, we obtain $\Projuivi M=\Projuivi (\Projuivi M)$ for any matrix $M \in \R^{n \times n}.$ In other words, $\Projuivi$ is a projection operator whenever one of the following is satisfied:
\begin{enumerate}[(i)]
\item $S$ is full column rank,
\item $T_j$ is full row rank for all $j \in \{1, \dots, m\},$
\item $T_1=T_2=\dots=T_m.$
\end{enumerate}

\noindent In \cite[Proposition 4.1]{hare2023hessianpublished}, if one of the previous three conditions is satisfied, it is shown that
\begin{align}\label{eq:doesnotaffect}
\Projuivi \shessti&=\shessti \quad \text{and} \quad \Projuivi \nabla_c^2 f(x^0;S,T_{1:m})=\nabla_c^2 f(x^0;S,T_{1:m}).
\end{align}

\noindent We use the following definition to describe the order at which an approximation technique converges.


\begin{definition} \cite[Definition 1.19]{Burden2016} \label{def:bigoh}
    Let $f:\R_+\to \R^{n \times p}$  and $g:\R_+\to \R.$ Suppose that $\lim_{\Delta \to 0} g(\Delta)=0$ and $\lim_{\Delta \to 0} f(\Delta)=L \in \R^{n \times p}.$ If there exists  a scalar $\kappa \geq 0$ with $$\Vert f(\Delta)-L \Vert \leq \kappa \, g(\Delta) \quad \text{for sufficiently small $\Delta$},$$
   then we say  $f(\Delta)$ is  $O(g(\Delta)),$ or  $f(\Delta)$ is a  $O(g(\Delta))$ accurate approximation of $L.$
\end{definition}

In this paper, $g(\Delta)$ takes the form $g(\Delta)=\Delta^N,$ where $N \in \N.$  
If Definition \ref{def:bigoh} is satisfied, we will  say that $f(\Delta)$ is an \emph{order-$N$ accurate approximation of $L$} where $N$ is the greatest positive integer satisfying  Definition \ref{def:bigoh}. We decided  to use the above definition rather than the definition of order-$N$ accuracy  as formulated in \cite[Definition 2.6]{hare2023hessianpublished}.  The above definition allows $L$ to be a square matrix or a  column vector which is  not possible in \cite[Definition 2.6]{hare2023hessianpublished}. We will see that $L$ takes the form of a column vector in Section \ref{sec:hvp}.

The following error bounds for the GSH and the GCSH were introduced in \cite{hare2023hessianpublished}. These error bounds  will be used  to  develop error bounds  for proper subsets of the entries of a Hessian. In the following theorems and the remainder of this paper, we use the notation
\begin{align*}
\Delta_u&=\max \{\Delta_S, \Delta_{T_1}, \dots, \Delta_{T_m}\}, \\ 
\Delta_l&=\min \{\Delta_S, \Delta_{T_1}, \dots, \Delta_{T_m}\}, \\ 
\widehat{T}&=\widehat T_j \quad \text{such that} \quad \left \Vert \widehat T_j^\dagger \right \Vert \quad \text{is maximal,} \quad j \in \{1, \dots, m\},  \\
k&= \max \{k_1, \dots, k_m\}. 
\end{align*}

\begin{theorem}[Error bounds for the GSH]\emph{\cite[Theorem 4.2]{hare2023hessianpublished}}\label{thm:mainprop}
Let $f:\dom f \subseteq \R^n\to\R$ be $\mathcal{C}^{3}$ on $B_n(x^0;\overline{\Delta})$ where $x^0 \in \dom f$ is the point of interest and $\overline{\Delta}>0$. Denote by $L_{\nabla^2 f}\geq 0$ the Lipschitz constant of $\nabla^2 f$ on $\overline{B}_n(x^0;\overline{\Delta})$. Let $S=\begin{bmatrix} s^1&s^2&\cdots&s^m\end{bmatrix} \in \R^{n \times m}$ and $T_j=\begin{bmatrix} t^1_j&t_j^2&\cdots&t_j^{k_j}\end{bmatrix} \in \R^{n \times k_j}$  for all $j \in \{1, \dots, m\}.$ Assume that $B_n(x^0;\Delta_{T_j})\subset B_n(x^0;\overline{\Delta})$ and $B_n(x^0+s^j;\Delta_{T_j})\subset B_n(x^0;\overline{\Delta})$  for all $j \in \{1, \dots m\}.$ Then the following hold.
 
 \begin{enumerate}[(i)]

 \item If $S$ is full column rank {\bf or}  $T_j$ is full row rank for all $j \in \{1, \dots, m\}$, then  
\begin{align}
\| \Projuivi \shessti- \Projuivi 
\nabla^2 f(x^0)\|&= \|  \shessti-  \Projuivi 
\nabla^2 f(x^0)\| \notag \\ &\leq 4m\sqrt{k} L_{\nabla^2 f}   \Vert (\widehat{S}^\top )^\dagger  \Vert \Vert \widehat{T}^\dagger \Vert   \left ( \frac{\Delta_u}{\Delta_l} \right )^2  \Delta_u.  \label{eq:seconditem}
\end{align}
\item  If $T_1=T_2=\dots=T_m=\T$, then 
\begin{align}
&\| \Projuv \shesst - \Projuv 
\nabla^2 f(x^0)\|= \|  \shesst- \Projuv  
\nabla^2 f(x^0)\| \notag \\
&\leq  4 \sqrt{mk} L_{\nabla^2 f} \frac{\Delta_u}{\Delta_l}  \left \Vert(\widehat{S}^\top)^\dagger \right \Vert  \left \Vert  \widehat{\T}^\dagger \right \Vert \Delta_u. \label{eq:simplexHessAllTequals}
\end{align}
\end{enumerate}
\end{theorem}

 Theorem \ref{thm:mainprop} shows that the GSH is an order-1 accurate approximation of the full Hessian. The accuracy of the GSH can be improved by decreasing both the radii $\Delta_u$ and $\Delta_\ell$ at the same rate.
\begin{theorem}[Error bounds for the GCSH]\emph{\cite[Theorem 4.3]{hare2023hessianpublished}} \label{thm:ebcentered2}
Let $f:\dom f \subseteq \R^n\to\R$ be $\mathcal{C}^{4}$ on $B_n(x^0;\overline{\Delta})$ where $x^0 \in \dom f$ is the point of interest and $\overline{\Delta}>0$. Denote by $L_{\nabla^3 f}$ the Lipschitz constant of $\nabla^3 f$ on $\overline{B}_n(x^0;\overline{\Delta}).$ Let  $S=\begin{bmatrix} s^1&s^2&\cdots&s^m\end{bmatrix} \in \R^{n \times m},$ $T_j=\begin{bmatrix} t_j^1&t_j^2&\cdots&t_j^{k_j} \end{bmatrix} \in \R^{n \times k_j}$ with the ball $B_n(x^0 + s^j;\Delta_{T_j})\subset B_n(x^0;\overline{\Delta})$  for all $j \in \{1, \dots, m\}$. Then the following hold.
\begin{enumerate}[(i)]

\item If $S$ is full column rank {\bf or}  $T_j$ is full row rank for all $j \in \{1, \dots, m\}$, then 
\begin{align}\label{eq:undertildecshesstwo}
\Vert \Projuivi \cshesstitwo - \Projuivi  
\nabla^2 f(x^0) \Vert&= \Vert  \cshesstitwo - \Projuivi 
\nabla^2 f(x^0) \Vert \notag \\
 &\leq 2m  \sqrt{k}L_{\nabla^3 f} \left (\frac{\Delta_u}{\Delta_l} \right )^2  \left \Vert  (\widehat{S}^\top)^\dagger  \right \Vert \left \Vert \left (\widehat{T} \right )^\dagger \right \Vert \Delta^2_u.
\end{align}
 \item If  $T_1=T_2=\dots=T_m=\T \in \R^{n \times k},$ then
\begin{align}
\left \|\Projuv \cshesst- \Projuv \nabla^2 f(x^0) \right \|&=\left \|  \cshesst- \Projuv \nabla^2 f(x^0) \right \| \notag \\
&\leq 2  \sqrt{mk}L_{\nabla^3 f} \frac{\Delta_u}{\Delta_l} \left \Vert  (\widehat{S}^\top)^\dagger  \right \Vert \left \Vert \left (\widehat{\T} \right )^\dagger \right \Vert \Delta^2_u. \label{eq:tildecentered2ebTiequaltwo} 
\end{align}
\end{enumerate}
\end{theorem}

 Theorem \ref{thm:ebcentered2} shows that the GCSH is an order-2 accurate approximation of the full Hessian. It provides a higher order of accuracy but it requires more function evaluations than the GSH when  $n \geq 2.$  Indeed, using a \emph{minimal poised set}, the GSH requires $(n+1)(n+2)/2$ function evaluations versus $n^2+n+1$ for the GCSH \cite[Section 5]{hare2023hessianpublished}. We are now ready to introduce the main results of this paper.

\section{Approximating a proper subset of the entries of the Hessian} \label{sec:partialHessian}

In this section,  we provide details on how  to choose the matrices of directions $S$ and $T_j$ when we are interested in a proper subset of the entries of the Hessian. In particular, we investigate how to approximate the diagonal entries, the off-diagonal entries and a column (row)  of the Hessian. The number of  function evaluations required is discussed and an error bound is provided in each case. Besides,   the relation between the \emph{centered simplex Hessian diagonal} (CSHD) introduced in \cite{jarry2022approximating} and the GCSH is discussed. The main contribution of this section  is to show how the GSH or the GCSH can be used to approximate a proper subset of the entries of a Hessian  and to provide error bounds showing the error is controllable in each case.  This shows that the GSH is a versatile tool that can either approximate a full Hessian or a proper subset of the entries of a Hessian. We begin by  presenting  results on how to approximate some or all  diagonal entries of a Hessian.

\subsection{Approximating the diagonal entries of the Hessian} \label{subsec:diagonalHessian}

An explicit formula to compute all the diagonal entries of the Hessian, which is well-defined  regardless of the number of sample points utilized, is discussed in \cite{jarry2022approximating}. The CSHD is an approximation technique that provides an order-2 accurate approximation of the diagonal entries of the Hessian.  We begin by showing that the CSHD  is a specific case of the GCSH when the appropriate matrices of directions $S$ and $T_j$ are employed. First, recall the definitions of the \emph{Hadamard product} and  the CSHD.

\begin{definition}\cite{Horn1990}
Let $A \in \R^{n \times m}$ and $B \in \R^{n \times m}.$ The Hadamard product of $A$ and $B$, denoted $A \odot B$ is the component-wise product. That is $[A \odot B]_{i,j}=A_{i,j}  B_{i,j}$ for all $i \in \{1, \dots n\}$ and $j \in \{1, \dots, m\}.$ 
\end{definition}

\begin{definition}[Centered simplex Hessian diagonal] \cite{jarry2022approximating} \label{def:cshd}
Let $f:\dom f \subseteq \R^n\to\R,$ $x^0 \in \dom f$ be the point of interest, $S=\begin{bmatrix} s^1&s^2&\cdots&s^m \end{bmatrix} \in \R^{n \times m}$ and $W=\begin{bmatrix} s^1 \odot s^1&\cdots& s^m \odot s^m  \end{bmatrix} \in \R^{n \times m}.$  Assume that $x^0 \oplus (\pm S)  \subset \dom f.$ 
 The \emph{centered simplex Hessian diagonal  of $f$ at $x^0$ over $S$},  denoted by $\cshd$ is  a vector in $\R^n$ given by 
\begin{align*}
    \cshd&=(W^\top)^\dagger \epsfS, \quad \text{where} \quad \epsfS=\begin{bmatrix} f(x^0+s^1)+f(x^0-s^1)-2f(x^0)\\ \vdots \\ f(x^0+s^m)+f(x^0-s^m)-2f(x^0)\end{bmatrix} \in \R^m.
\end{align*}
\end{definition}

 A specific type of matrix is involved in this section, called \emph{partial diagonal matrix}. The definition follows.
\begin{definition}[Partial diagonal matrix] \label{def:partial}
Let $M \in \R^{n \times m}, m \leq n.$ We say that $M$ is a partial diagonal matrix if   there exists a diagonal matrix $D \in \R^{n \times n}$ such that for each column $Me^j_m, j \in \{1, \dots m\},$ there exists a unique distinct index $i \in \{1, \dots, n\}$ that yields $Me^j_m=De^i_n.$
\end{definition}

 In other words, a partial diagonal matrix is a subset of the columns of a  single diagonal matrix. Note that the columns in a partial diagonal matrix do not need to be ``ordered''. For example, the matrix $$M=\begin{bmatrix} 1&0\\0&0\\0&2\end{bmatrix} \quad \text{and} \quad \ddot{M}=\begin{bmatrix} 0&1\\0&0\\2&0 \end{bmatrix}$$ are  partial diagonal matrices, but $$\widetilde{M}=\begin{bmatrix} 1&0&3\\0&0&0\\0&2&0\end{bmatrix}, \quad \overline{M}=\begin{bmatrix} 1&1\\0&0\\0&0\end{bmatrix}$$ are not  partial diagonal matrices. 
Note that a  partial diagonal matrix is full column rank if and only if  it does not contain a column  equal to the zero vector in $\R^n.$ The following lemma provides details about  the Moore--Penrose pseudoinverse of a partial diagonal matrix with full column rank.

\begin{lemma} \label{lem:partialdiag}
Let $S=\begin{bmatrix} s^1&s^2&\cdots&s^m \end{bmatrix} \in \R^{n \times m}$ where $m \leq n$  be a partial  diagonal matrix with full column rank. Then $$S^\dagger=\begin{bmatrix} (s^1)^\dagger \\(s^2)^\dagger \\\vdots \\(s^m)^\dagger\end{bmatrix} .$$
\end{lemma}
\begin{proof}
Let $u_j$ be the index in $\{1, \dots, m\}$ of the only non-zero entry in column $s^j.$ Since $S$ is full column rank, using \eqref{eq:fullcolumnrank}, we have 
$$
S^\dagger=(S^\top S)^{-1}S^\top=\left (\Diag \begin{bmatrix} (s^1_{u_1})^2&\cdots & (s^m_{u_m})^2 \end{bmatrix} \right )^{-1}S^\top=\Diag \begin{bmatrix} \frac{1}{(s^1_{u_1})^2}&\cdots & \frac{1}{(s^m_{u_m})^2} \end{bmatrix} S^\top=\begin{bmatrix} \frac{1}{s^1_{u_1}}(e^{u_1})^\top\\ \frac{1}{s^2_{u_2}}(e^{u_2})^\top\\ \vdots \\ \frac{1}{s^m_{u_m}}(e^{u_m})^\top\end{bmatrix}.
$$

\noindent Since $S$ is full column rank, using \eqref{eq:fullcolumnrank} we find  $(s^j)^\dagger=\frac{1}{s^j_{u_j}}(e^{u_j})^\top$ for all $j \in \{1, \dots, m\}$ and the result follows. \qed
\end{proof}

  The following theorem provides a sufficient condition for the GCSH to return the same approximation of the diagonal entries of the Hessian as the CSHD.

\begin{theorem}\label{thm:shequalcshd}
Let $f:\dom f \subseteq \R^n \to \R, x^0 \in \dom f$ be the point of interest,  $S=\begin{bmatrix} s^1&s^2&\cdots&s^m \end{bmatrix} \in \R^{n \times m}$ and $T_{j}=-s^j \in \R^n$ for all  $j \in \{1, \dots, m\}.$ Let $z \in \R^n$ be a vector containing the $n$ diagonal entries of $\cshtwo  f(x^0;S,T_{1:m})$. That is $z_i=\left [\cshtwo f(x^0;S,T_{1:m})\right ]_{i,i}$ for all $i \in \{1, \dots, n\}.$    If $S$ is a  partial diagonal matrix with full column rank, then $z=\cshd.$
\end{theorem}
\begin{proof}
Let $A=\begin{bmatrix} S&-S \end{bmatrix} \in \R^{n \times 2m}$ and $T_{m+j}=s^j$  for $j \in \{1, \dots, m\}.$ We have
\begin{align*}
\cshtwo f(x^0;S,T_{1:m})&=\nabla_s^2 f(x^0;A,T_{1:2m}) \quad \quad \text{(by Proposition 5.9 in \cite{jarry2023numerical})} \\
&=(A^\top)^\dagger \begin{bmatrix} \left ( \nabla_s f(x^0+s^1;-s^1)-\nabla_s f(x^0;-s^1) \right )^\top \\ \vdots \\ \left ( \nabla_s f(x^0+s^m;-s^m)-\nabla_s f(x^0;-s^m) \right )^\top \\ \left ( \nabla_s f(x^0-s^1;s^1)-\nabla_s f(x^0;s^1) \right )^\top \\ \vdots \\ \left ( \nabla_s f(x^0-s^m;s^m)-\nabla_s f(x^0;s^m) \right )^\top \end{bmatrix}. \\
\end{align*}
Since $(A^\top)^\dagger=\frac{1}{2} \begin{bmatrix} (S^\top)^\dagger &-(S^\top)^\dagger\end{bmatrix},$ and  expanding each row  of the form  $$\left ( \nabla_s f(x^0 \pm s^j; \mp s^j)-\nabla_s f(x^0; \mp s^j) \right )^\top,$$ we obtain
\begin{align*}
\cshtwo f(x^0;S,T_{1:m})&=\frac{1}{2} \begin{bmatrix} (S^\top)^\dagger &\, -(S^\top)^\dagger\end{bmatrix} \begin{bmatrix} (-s^1)^\dagger \left ( -f(x^0+s^1)-f(x^0-s^1)+2f(x^0) \right )\\ \vdots \\ (-s^m)^\dagger \left ( -f(x^0+s^m)-f(x^0-s^m)+2f(x^0)\right )\\ (s^1)^\dagger \left ( 2f(x^0)-f(x^0-s^1)-f(x^0+s^1)\right )\\ \vdots \\ (s^m)^\dagger \left ( 2f(x^0)-f(x^0-s^m)-f(x^0+s^m)\right ) \end{bmatrix} \\
&=(S^\dagger)^\top \begin{bmatrix} (s^1)^\dagger \left ( f(x^0-s^1)+f(x^0+s^1)-2f(x^0) \right )\\ \vdots \\ (s^m)^\dagger \left ( f(x^0-s^m)+f(x^0+s^m)-2f(x^0)\right ) \end{bmatrix} \\
&=\begin{bmatrix} ((s^1)^\dagger)^\top&\cdots&((s^m)^\dagger)^\top \end{bmatrix} \begin{bmatrix} (s^1)^\dagger \left ( f(x^0-s^1)+f(x^0+s^1)-2f(x^0) \right )\\ \vdots \\ (s^m)^\dagger \left ( f(x^0-s^m)+f(x^0+s^m)-2f(x^0)\right ) \end{bmatrix}
\end{align*}
by Lemma \ref{lem:partialdiag}. Let $z \in \R^n$ be the vector containing the $n$ diagonal entries  of the previous equation. Then 
\begin{align*}
z&=\begin{bmatrix} ((s^1)^\top)^\dagger \odot ((s^1)^\top)^\dagger &\cdots&((s^m)^\top)^\dagger \odot ((s^m)^\top)^\dagger\end{bmatrix}  \epsfS \\
&=(W^\top)^\dagger \epsfS=\cshd.
\end{align*} \qed
\end{proof} 

By defining the sets $T_{j}$   as in Theorem \ref{thm:shequalcshd}, the CSHD and the GCSH  use the  same set of sample points. However, if $S$ is not a partial diagonal matrix with full column rank, then  the vector $z$ containing the diagonal entries of the GCSH is not necessarily equal to the CSHD.  Moreover, the  GCSH  is not necessarily a diagonal matrix.  The following two examples illustrate these claims.

\begin{example}
Let $$S=\begin{bmatrix} s^1&s^2&s^3 \end{bmatrix}=\begin{bmatrix} 0.1&0&0\\0&0.1&0.2\\0&0&0 \end{bmatrix}.$$ Let $T_{j}=-s^j$ for all $j \in \{1, 2, 3\}$.
Let  $f(x)=-2x_1^4+x_2^4+10x_3^4$ and $x^0=\begin{bmatrix}2&-2&5\end{bmatrix}^\top.$  Note that $$\nabla^2 f(x^0)=\Diag [-96~48~3000].$$ The GCSH is  $$\nabla^2_c 
 f(x^0;S,T_{1:3})= \Diag [-96.04~48.068~0],$$  and the CSHD is  $$\cshd=\begin{bmatrix} -96.04&48.0765&0 \end{bmatrix}^\top.$$
\end{example}

The next example shows that the GCSH is not necessarily a diagonal matrix, even when we use the same set of sample points.
\begin{example}
Let  $$S=\begin{bmatrix} s^1&s^2\end{bmatrix}=\begin{bmatrix} 0.1&0.1\\0&0.1\\0&0 \end{bmatrix}.$$ Let $T_{j}=-s^j$ for all $j \in \{1, 2\}.$ Consider the same function and point of interest as in the previous example. That is  $f(x)=-2x_1^4+x_2^4+10x_3^4$ and $x^0=\begin{bmatrix}2&-2&5\end{bmatrix}^\top.$   Then the GCSH is  $$\nabla^2_c 
 f(x^0;S,T_{1:2})= \begin{bmatrix} -96.04&0&0\\72.03&-24.01&0\\0&0&0\end{bmatrix},$$ and  the CSHD is $$\cshd=\begin{bmatrix} -96.04&48.02&0 \end{bmatrix}^\top.$$
\end{example}

To sum up, the CSHD or the GCSH can both be used to approximate the diagonal entries of a Hessian.  It returns the same values for the diagonal entries  whenever $S$ is a partial diagonal matrix with full column rank.  In general, both techniques do not necessarily return the same approximation of  the diagonal entries of a Hessian. 

The following proposition provides a  result concerning the projection of a matrix over $S$ and $T_{1:m}.$

\begin{proposition} \label{lem:projHessianequalprojD}
Let $M \in  \R^{n \times n}.$ Let   $S=\begin{bmatrix} s^1&\cdots&s^m \end{bmatrix} \in \R^{n \times m}$ and   let $T_{j}=-s^j$ for all $j \in \{1, \dots, m\}.$  If $S$ is a partial diagonal matrix with full column rank, then
\begin{align*}
\Projuivi M&= \Projuivi \Diag \begin{bmatrix} M_{1,1}&\cdots&M_{n,n} \end{bmatrix}.
\end{align*}
Moreover,  if $(e^i)^\top S\neq \zero_m^\top$ for some $i \in \{1, \dots, n\},$ then  
\begin{align*}
     \left [\Projuivi M \right ]_{i,i}&=M_{i,i}. 
\end{align*}
If  $(e^i)^\top S= \zero_m^\top$ for some $i \in \{1, \dots, n\},$ then $$\left [\Projuivi M \right ]_{i,i}=0.$$ 
\end{proposition}
\begin{proof}
We have
\begin{align*}
\sum_{j=1}^m (S^\top)^\dagger e^j(e^j)^\top S^\top M T_jT_j^\dagger&=\sum_{j=1}^m ((s^j)^\top)^\dagger (s^j)^\top M (-s^j)(-s^j)^\dagger \\
&=\sum_{j=1}^m e^{u_j} (e^{u_j})^\top M e^{u_j} (e^{u_j})^\top
\end{align*}
where $u_j$ represents the  index of the only non-zero entry in $s^j,$  $u_j \in \{1, \dots n\},$ and $j \in \{1, \dots, m\}.$ From the definition of a  partial diagonal matrix, we know that $u_j \neq u_{\bar{j}}$  whenever $j \neq \bar{j},$ $j$ and $\bar{j}$ in $\{1, \dots, m\}$. Noticing that $e^{u_j}(e^{u_j})^\top=\Diag(e^{u_j}),$ we get
\begin{align*}
\sum_{j=1}^m (S^\top)^\dagger e^j(e^j)^\top S^\top M T_jT_j^\dagger&=\sum_{j=1}^m \Diag(e^{u_j}) M \Diag(e^{u_j})\\
&=\sum_{j=1}^m \Diag(e^{u_j})\cdot M_{u_j,u_j}=\Projuivi \Diag[M_{1,1}~\cdots~M_{n,n}]. 
\end{align*} 
The rest of the proof follows immediately  from the fact that $m \leq n,$ and $u_j \neq u_{\bar{j}}$ whenever $j \neq \bar{j},$ $j$ and $\bar{j}$ in $\{1, \dots, m\}.$ \qed
\end{proof} 

The notation $D$ is now used to represent the diagonal matrix in $\R^{n \times n}$ containing the diagonal entries of the Hessian $\nabla^2 f(x^0).$ That is $D_{i,i}=[\nabla^2 f(x^0)]_{i,i}$ for all $i \in \{1, \dots, n\}.$ 
If $S$ is a  diagonal matrix with full column rank and $T_j=-s^j$ for all $j \in \{1, \dots, n\},$ it follows from Proposition \ref{lem:projHessianequalprojD} that $$\Projuivi \nabla^2 f(x^0)=\Projuivi D.$$  In other words, the projection of the full true Hessian  is a diagonal matrix that keeps intact all diagonal entries of the true Hessian. In the case where $S$ is  a non-square partial diagonal matrix, then  it   makes the $(i,i)$ diagonal entry of the true Hessian equal to zero if $S$ does not contain a multiple of the identity column $e^i.$ Also, since $S$ is full column rank, it follows from Propositions \ref{lem:projHessianequalprojD} and Equation \eqref{eq:doesnotaffect}  that $\shessti$ and $\cshesstitwo$ are diagonal matrices.

The next theorem presents an error bound when the GCSH is used to approximate some, or all  diagonal entries of the true Hessian. 

 \begin{corollary}[Error bound for the diagonal entries of the Hessian] \label{thm:ebdiagentries}
Let $f:\dom f\subseteq \R^n \to \R$ be $\mathcal{C}^4$ on an open domain containing $\overline{B}_n(x^0;\Delta_S)$ where $x^0 \in \dom f$ is the point of interest and $\Delta_S>0$ is the radius of $S=\begin{bmatrix} s^1&\cdots s^{m}\end{bmatrix} \in \R^{n \times m}.$  Let $T_{j}=-s^j$ for all $j \in \{1, \dots, m\}.$  Denote by $L_{\nabla^3 f} \geq 0$ the Lipschitz constant of $\nabla^3 f$ on $\overline{B}_n(x^0;\Delta_S)$. If $S$ is a partial diagonal matrix with full column rank, then  
\begin{align}
   &\left \Vert  \Projuivi \cshesstitwo -\Projuivi \nabla^2 f(x^0) \right \Vert= \left \Vert \cshesstitwo-\Projuivi D \right \Vert \leq   \frac{1}{12}L_{\nabla^3 f} \Delta_S^2.\label{eq:EBdiag}
\end{align}
\end{corollary}

\begin{proof}
By Equation \eqref{eq:doesnotaffect} and Proposition \ref{lem:projHessianequalprojD},  we get the equality. 
To make notation more compact, let $\varepsilon=\epsfS \in \R^m.$ We have
\begin{align*}
\left \Vert  \cshesstitwo- \Projuivi D \right \Vert&= \scalemath{1}{\left \Vert \sum_{i=1}^m ((s^i)^\top)^\dagger (s^i)^\top (S^\top)^\dagger \Diag(\varepsilon) S^\dagger (-s^i)(-s^i)^\dagger-\sum_{i=1}^m ((s^i)^\top)^\dagger (s^i)^\top D s^i (s^i)^\dagger   \right \Vert} \\
& \leq \max_{i=1, \dots, m} \left ( \Vert ((s^i)^\top)^\dagger \Vert \Vert (s^i)^\dagger \Vert \left \vert (s^i)^\top (S^\top)^\dagger \Diag(\epsilon) S^\dagger s^i-(s^i)^\top Ds^i \right \vert  \right )\\
&=\max_{j=1, \dots, m} \left ( \frac{1}{\Vert s^j \Vert^2 } \left \vert \varepsilon_j -(s^j)^\top Ds^j \right \vert  \right ).
\end{align*}
 By Taylor's Theorem, using a similar process as in the proof in \cite[Theorem 3.3]{jarry2022approximating}, we  obtain 
\begin{align*}
\left \Vert  \cshesstitwo- \Projuivi D \right \Vert &\leq \max_{j=1, \dots, m} \left ( \frac{1}{\Vert s^j \Vert^2 } \frac{1}{12}L_{\nabla^3 f} \Vert s^j \Vert^4 \right ) \\
&= \max_{j=1, \dots, m} \left ( \frac{1}{12}L_{\nabla^3 f} \Vert s^j \Vert^2 \right )\\
&\leq \frac{1}{12}L_{\nabla^3 f} \Delta_S^2. 
\end{align*} \qed
\end{proof}

By defining $S$ and $T_j$ as in the previous corollary,  note that the general error bound proposed for the GCSH in Theorem \ref{thm:ebcentered2}$(ii)$  is also valid. The previous proof utilized properties of partial diagonal matrices to obtain a tighter  error bound than the one proposed in Theorem \ref{thm:ebcentered2}$(ii)$.

The previous corollary shows how to obtain an order-2 accurate approximation of some, or all diagonal entries of the Hessian.  This  requires  $2n+1$ function evaluations when $S$ is square.  
\begin{figure}[h]
\centering
\includegraphics[scale=0.35]{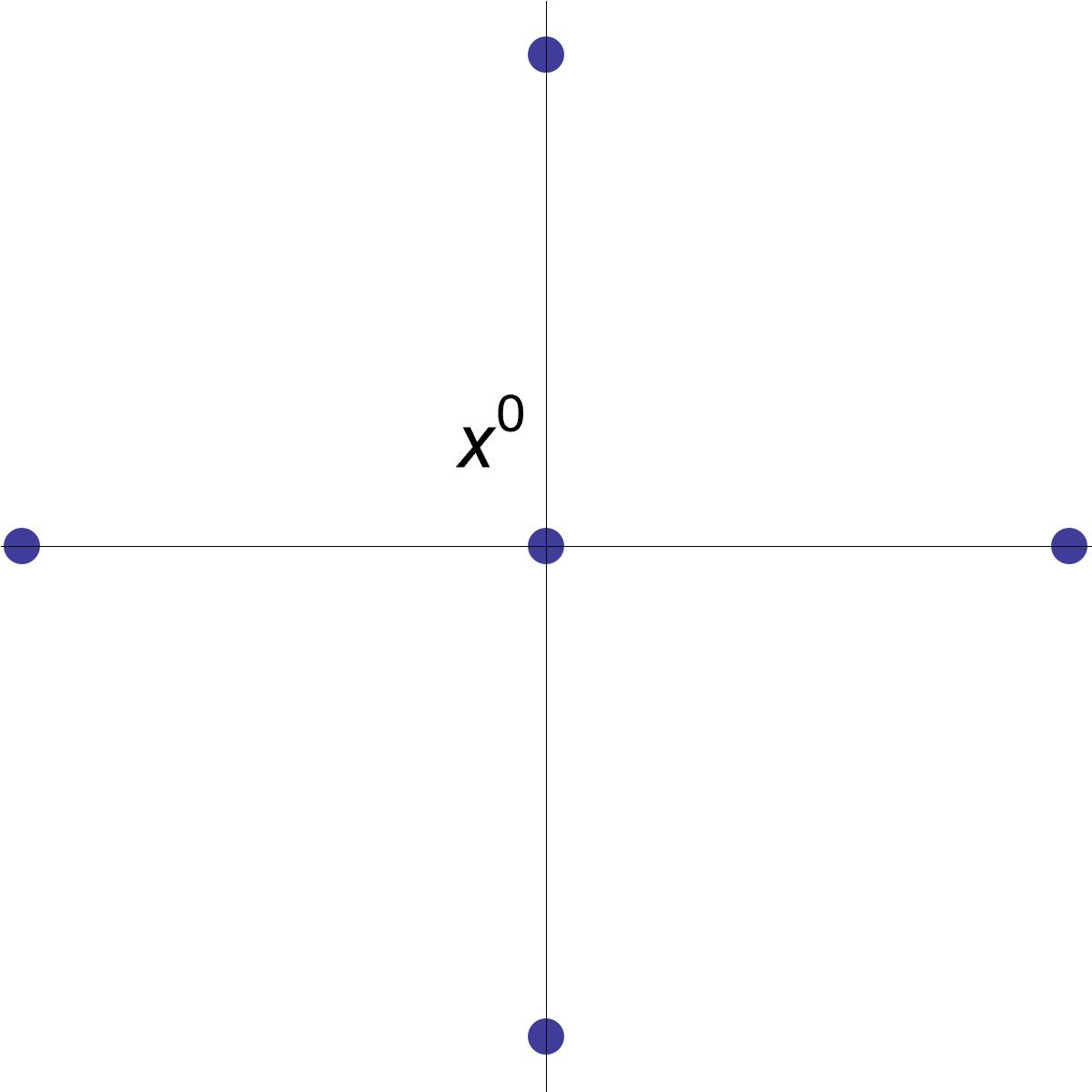}
\caption{Sample set created while computing the GCSH with $S=\Id_2, T_1=-e^1_2, T_2=-e^2_2$ in Corollary \ref{thm:ebdiagentries}}
\end{figure} 

\noindent If we are interested  in approximating only one diagonal entry of a Hessian $\nabla^2 f(x^0)$, say $\left [\nabla^2 f(x^0) \right ]_{i,i},$  then the computational cost is  three function evaluations. In this case, we   can choose $S=he^i$ and $T_1=-he^i.$  Each additional   diagonal entry can  be obtained for  two  more function evaluations.

Other  matrices of directions $S$ and $T_{j}$ may be used to obtain an approximation of all diagonal entries of a Hessian. For instance,  the following matrices  can be used:
$$S=\begin{bmatrix} s^1&\cdots&s^n \end{bmatrix} =h\Id, \quad T_{j}=s^j, \quad \text{for all} \, j \in \{1, \dots, n\}, h \neq 0.$$
In this case, $S$ is diagonal with full column rank and it follows from Proposition \ref{lem:projHessianequalprojD} that  $\Projuivi \nabla^2 f(x^0)=\Projuivi D=D$ where $D \in \R^{n \times n}$ is the diagonal matrix such that $D_{i,i}=[\nabla^2 f(x^0)]_{i,i}$ for all $i \in \{1, \dots, n\}.$  By Theorem \ref{thm:mainprop}$(ii),$ this choice of matrices provides an order-1 accurate approximation of all diagonal entries of the Hessian. The computation of $\nabla_s^2 f(x^0;S,T_{1:n})$  requires $2n+1$ function evaluations.  Hence, it is preferable to choose the matrices of directions $S$ and $T_{j}$ as in Corollary \ref{thm:ebdiagentries}, since  it provides a greater order of accuracy for the same number of function evaluations.

In the next section, we investigate the approximation of some, or all off-diagonal entries of the Hessian.

\subsection{Approximating the off-diagonal entries of the Hessian} \label{subsec:offDiagonalHessian}
 In this section, how to approximate some, or all off-diagonal entries of the Hessian is examined. First, recall that  the Hessian $\nabla^2 f(x^0)$ is symmetric whenever $f \in \mathcal{C}^2$. Therefore,  it is sufficient to  consider  the off-diagonal entries     $[\nabla^2 f(x^0)]_{i,j}$ such that $i<j.$  It is possible to approximate some, or all off-diagonal entries of the Hessian  by setting the matrices of directions $S$ and $T_j$ in the following way. Define
\begin{align}
&\widetilde{S} \in \R^{n \times n-1}:  \text{a  partial diagonal matrix with full column rank} \notag \\
&\quad \quad \quad \quad \quad \quad \, \text{such that  the $n$\textsuperscript{th} row  is   equal to $\zero_{n-1}^\top,$} \notag \\
&S=\begin{bmatrix} s^1&\cdots&s^m\end{bmatrix} \in \R^{n \times m}: \text{a  non-empty subset of the columns of $\widetilde{S},$}  \label{eq:Soffdiag}\\
&T=\begin{bmatrix} t^1&\cdots&t^n \end{bmatrix} \in \R^{n \times n}: \text{a diagonal matrix with full column rank,} \notag \\
&\widetilde{T}_j=\begin{bmatrix} t^{u_j+1}&\cdots&t^{n-1}&t^n\end{bmatrix} \in \R^{n \times n-u_j} \text{where $u_j$ represents the index} \notag \\
& \text{of the non-zero entry in $s^j,$} \, j \in \{1, \dots, m \}, \notag \\
&T_{j} \in \R^{n \times k_j}: \text{a subset of directions contained in $\widetilde{T}_j$ for all $j \in \{1, \dots, m\}.$} \label{eq:Tioffdiag}
\end{align}
In the next corollary, the matrix $U \in \R^{n \times n}$ denotes a strictly upper triangular matrix such that \begin{equation*} 
  U_{i,j}=\left\{
  \begin{aligned}
  \left [\nabla^2 f(x^0)\right ]_{i,j}&, && \text{if} \quad  1 \leq i<j \leq n, \\
    &0, && \text{otherwise}.
  \end{aligned}\right.
\end{equation*} 

\noindent Using a similar process to the one in Proposition \ref{lem:projHessianequalprojD}, it can be shown that $$\Projuivi \nabla^2 f(x^0)=\Projuivi U$$ and that the GSH (GCSH) is a strictly upper triangular matrix   
whenever the matrices of directions $S$ and  $T_j$ are defined as in \eqref{eq:Soffdiag} and \eqref{eq:Tioffdiag}.

The following  two error bounds  follow from   Theorem \ref{thm:mainprop}$(ii)$ and Theorem \ref{thm:ebcentered2}$(ii)$, respectively. 

\begin{corollary}[Error bound for the off-diagonal entries of the Hessian] \label{thm:eboffdiag}
Let $f:\dom f \subseteq \R^n\to\R$ be $\mathcal{C}^{4}$ on $B_n(x^0;\overline{\Delta})$ where $x^0 \in \dom f$ is the point of interest and $\overline{\Delta}>0$. Denote by  $L_{\nabla^2 f}\geq 0$ and $L_{\nabla^3 f} \geq 0$ the Lipschitz constant of $\nabla^2 f$ and $\nabla^3 f$  on $\overline{B}_n(x^0;\overline{\Delta})$ respectively.  Let $S=\begin{bmatrix} s^1&s^2&\cdots&s^m\end{bmatrix} \in \R^{n \times m}$ and $T_j \in \R^{n \times k_j}$  be defined as in \eqref{eq:Soffdiag} and \eqref{eq:Tioffdiag} respectively.  Assume that $B_n(x^0+s^j;\Delta_{T_j})\subset B_n(x^0;\overline{\Delta})$  for all $j \in \{1, \dots, m\}.$ Then
\begin{enumerate}[(i)]
\item 
\begin{align*}
\left \Vert \Projuivi \shessti- \Projuivi \nabla^2 f(x^0)  \right \Vert&=\Vert \shessti- \Projuivi U \Vert \\
&\leq   4m\sqrt{k} L_{\nabla^2 f} \left (\frac{\Delta_u}{\Delta_l} \right )^2 \Vert (\widehat{S}^\top)^\dagger \Vert \Vert \widehat{T}^\dagger \Vert \Delta_u,
 \end{align*}
 and
\item 
\begin{align*}
\left \Vert \Projuivi \cshesstitwo- \Projuivi \nabla^2 f(x^0)  \right \Vert&=\left \Vert \cshesstitwo- \Projuivi U \right \Vert  \\
&\leq  2m  \sqrt{k}L_{\nabla^3 f} \left (\frac{\Delta_u}{\Delta_l} \right )^2 \Vert (\widehat{S}^\top)^\dagger \Vert \Vert \widehat{T}^\dagger \Vert \Delta^2_u.
\end{align*}
\end{enumerate}
\end{corollary}

 A simple choice for $S$ and $T_j$  if \textbf{all} off-diagonal entries  are of interest is to set 
\begin{align}
S&=h\begin{bmatrix} e^1&\cdots&e^{n-1} \end{bmatrix}, \label{eq:Soffdiagex} \\
T_j&=h \begin{bmatrix} e^{j+1}&\cdots&e^n \end{bmatrix}, \quad \text{for all} \quad j \in \{1, \dots, n-1\} \label{eq:offdiag}
\end{align}
where $h \neq 0.$ In this case, the GSH is an order-1 accurate approximation of all off-diagonal entries of the Hessian. To compute this GSH, the function  must be evaluated at the points $x^0, x^0 \oplus S, x^0 \oplus T_j$ and $ x^0 +s^j \oplus T_j$ for all $j \in \{1, \dots n-1\}.$ Hence, the number of distinct function evaluations is 
\begin{equation*}
\scalemath{0.95}{1+(n-1)+(n-1)+\frac{(n-1)n}{2}-(n-2)=n+\frac{(n-1)n}{2}=\frac{n(n+1)+2}{2}.}
\end{equation*}
In the previous equation, we subtracted $(n-2)$ since $x^0 \oplus h \begin{bmatrix} e^2&\cdots e^{n-1} \end{bmatrix}$  appears in $x^0 \oplus T_j$ and $x^0 \oplus S.$ Note that this  number of function evaluations is smaller than  $(n+1)(n+2)/2$  whenever $n \geq 1,$ which is the number of function evaluations require to compute  a GSH  with a \emph{minimal poised set for GSH} \cite[Definition 5.2]{hare2023hessianpublished}. Therefore, if we are only interested in  the off-diagonal entries of a Hessian, it is preferable to set the matrices $S$ and $T_j$ as described  in this section, rather than using a minimal poised set for GSH.

\begin{figure}[ht]
\centering
\includegraphics[scale=0.35]{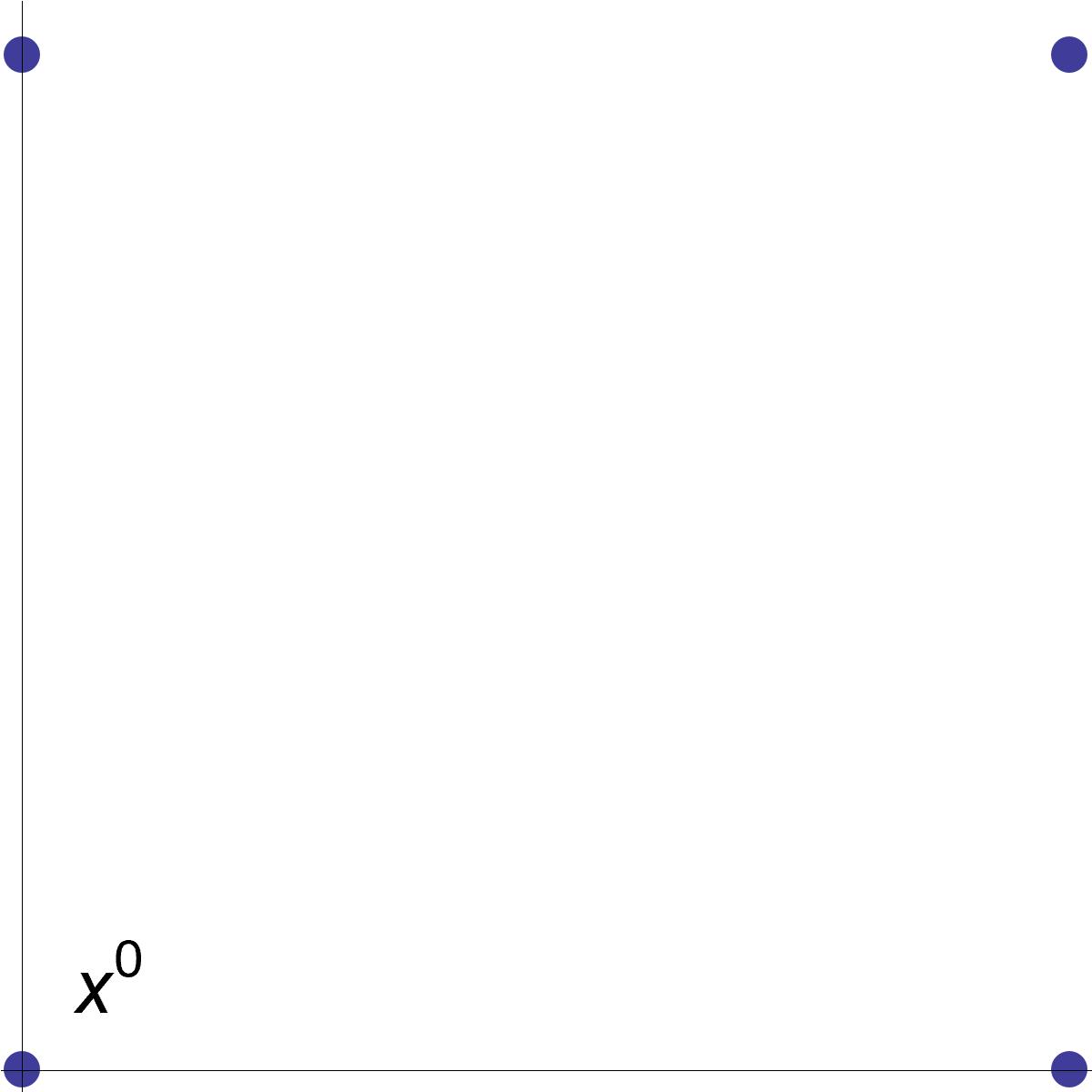}
\caption{Sample set  created while computing the GSH with  $S=e^1_2,T_1=e^2_2$ as described in Eqs. \eqref{eq:Soffdiagex} and \eqref{eq:offdiag}}
\end{figure} 
 
In the previous corollary, item $(ii)$ shows  that the GCSH  is an order-2 accurate approximation of all off-diagonal entries of the Hessian. In this case, the sample points used are $x^0, x^0 \oplus  (\pm S), x^0 \oplus (\pm T_j), x^0 \oplus  S \oplus T_j,$ and $x^0 \oplus (-S) \oplus (-T_j)$ for all $j \in \{1, \dots n-1\}.$  The number of distinct  function evaluations is 
\begin{align*}
1+ 2 \left ( \frac{n(n+1)}{2} \right )&=n^2+n+1.
\end{align*}
Notice that an approximation of the full Hessian can be  obtained with $n^2+n+1$  function evaluations by  taking $S \in \R^{n \times n}$ full rank and $\overline{T}=-S.$ We refer to this choice as a \emph{minimal poised set for GCSH}. Therefore, there is no advantage in terms of function evaluations to choose $S$ and $T_j$ as described in \eqref{eq:Soffdiag} and \eqref{eq:Tioffdiag} over a minimal poised set for GCSH.

\begin{figure}[h]
\centering
\begin{minipage}{.45\textwidth}
  \centering
  \includegraphics[width=.6\linewidth]{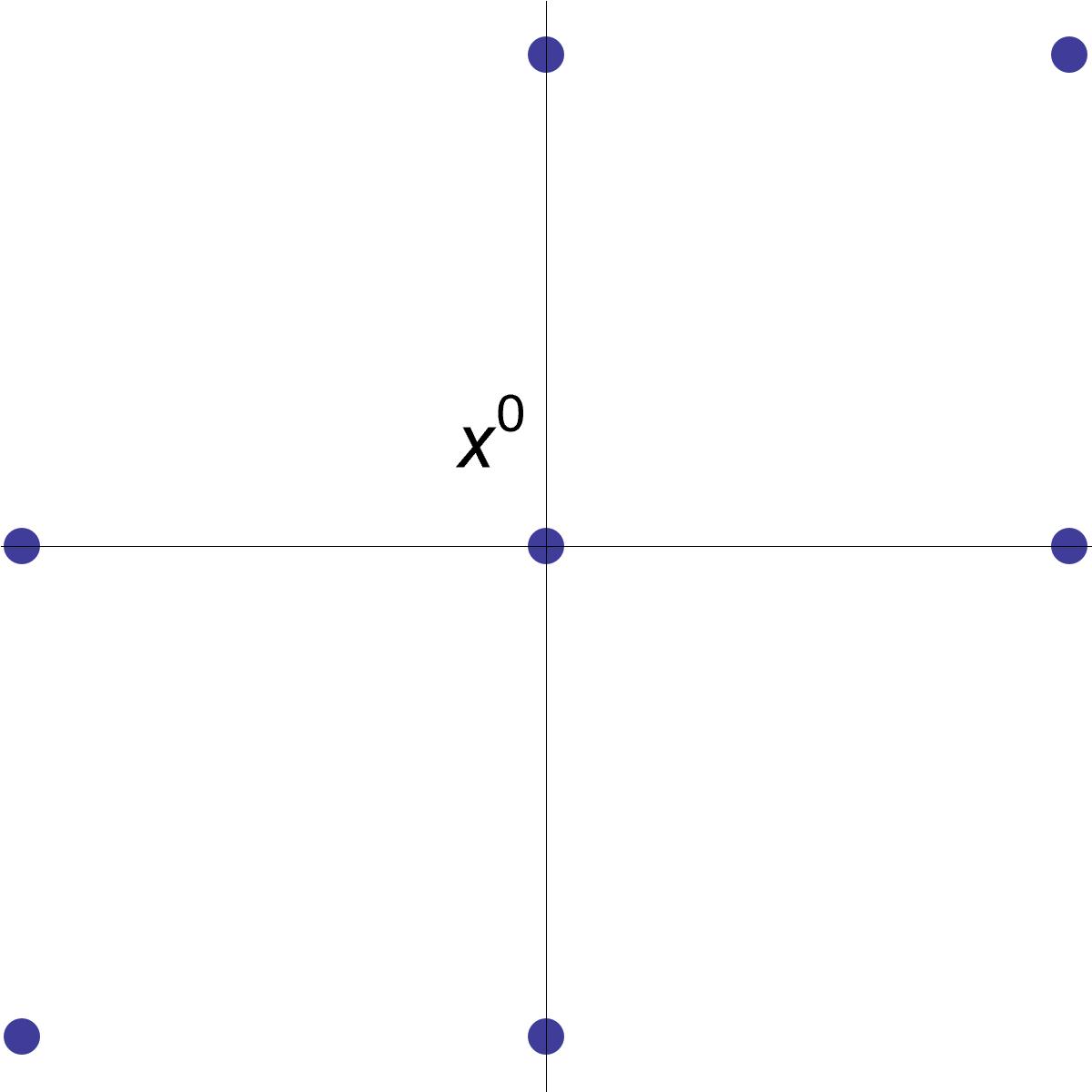}
  \caption{Sample set  created while computing the GCSH with $S=e^1_2$ and $T_1=e^2_2$ }
  \label{fig:test1}
\end{minipage}%
\hspace{0.3cm}
\begin{minipage}{.45\textwidth}
  \centering
  \includegraphics[width=.6\linewidth]{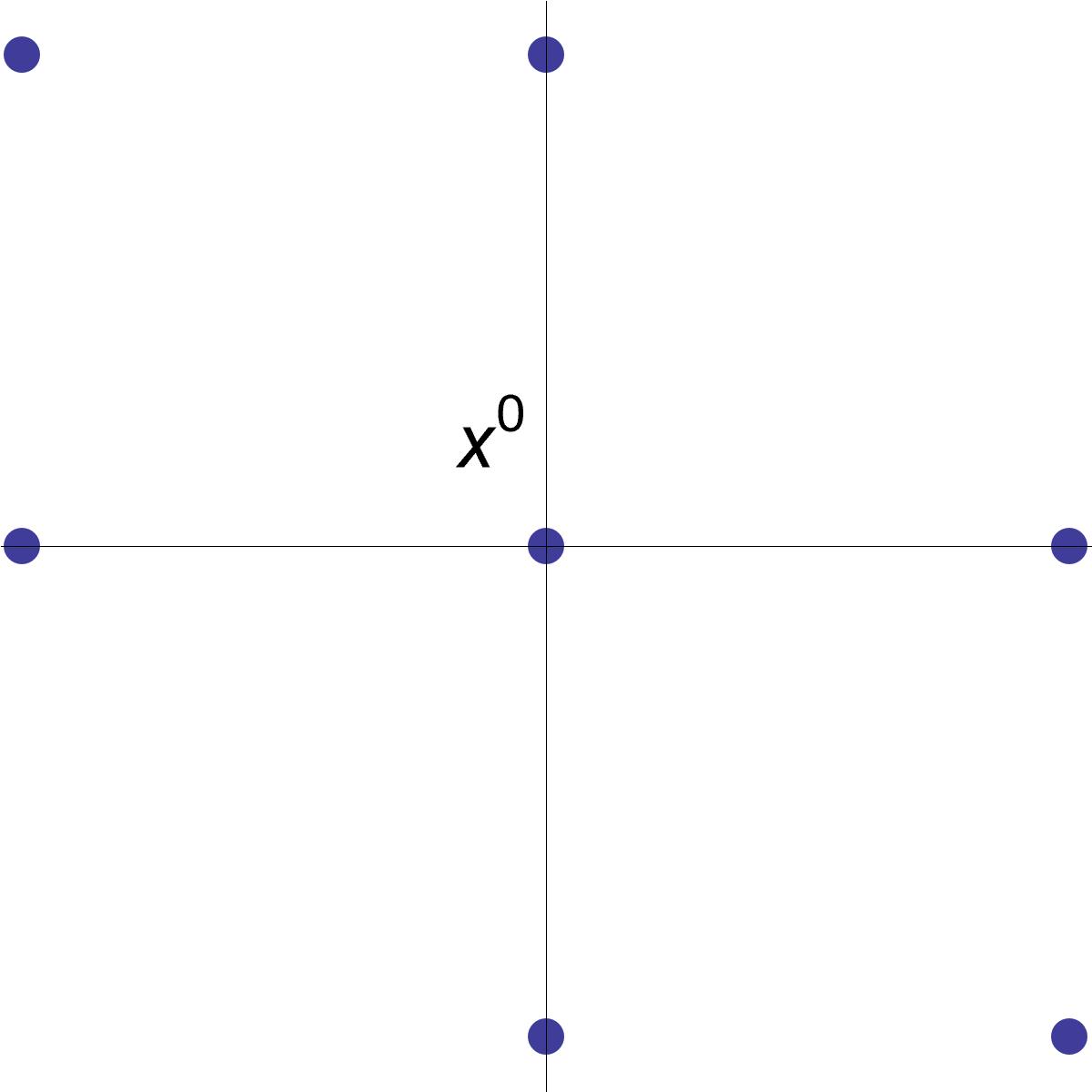}
  \caption{Sample set created while computing the GCSH with $S=\Id_2$ and $\T=-\Id_2$}
  \label{fig:test2}
\end{minipage}
\end{figure}

It does not seem possible to obtain an order-1 accurate approximation of all off-diagonal entries of a Hessian with fewer than $\frac{n(n+1)}{2}+1$ function evaluations, nor an order-2 accurate approximation with fewer than $n^2+n+1$ function evaluations. An obvious future research direction is to investigate this conjecture and mathematically prove or disprove it.
In the next section, we discuss how to approximate one row of the Hessian. 

\subsection{Approximating a row/column of the Hessian} \label{subsec:ColumnHessian}

In this section, we discuss how to approximate  some, or all entries of a row/column in the Hessian. Since the Hessian is symmetric, approximating a row also provides an approximation of the corresponding column.  

Let $M \in \R^{n \times n}.$ We denote by $R_i \in \R^{n \times n}$ the square matrix such that  $R_i=\Diag(e^i) M$  for all $i \in \{1, \dots, n\}$.  We begin  by introducing the following lemma.
\begin{lemma}
Let $M \in \R^{n \times n},$ $S=he^i \in \R^n$ where $h\neq 0,$  and $\T  \in \R^{n \times k}.$ Define $R_i=\Diag(e^i) M$  for all $i \in \{1, \dots, n\}$.  Then for all $i \in \{1, \dots, n\},$
$$\Projuv M= \Projuv R_i.$$
\end{lemma}
\begin{proof}
 We have
\begin{align*}
\Projuv M&=((he^i)^\top)^\dagger (he^i)^\top M \T\T^\dagger \\
&=e^i(e^i)^\top M \T\T^\dagger \\
&=R_i \T\T^\dagger \\
&=(e^i)(e^i)^\top R_i \T\T^\dagger \\
&= ((he^i)^\top)^\dagger (he^i)^\top R_i \T\T^\dagger=\Projuv R_i. 
\end{align*} \qed
\end{proof}

In words, the previous result says that the projection onto $S$ and $\T$ of a matrix $M$ is equal to the projection onto $S$ and $\T$ of row $i$  of this matrix  whenever $S=he^i.$

When $S$ and $\T$ are defined as in the previous lemma, $S$ is full column rank and it follows from Equation \eqref{eq:doesnotaffect}  that $\shesst=\Diag(e^i)\shesst$ and $\cshesst=\Diag(e^i)\cshesst.$ Moreover, the  projection of the Hessian is $$\Projuv \nabla^2 f(x^0)=\Projuv \Diag(e^i)\nabla^2 f(x^0)$$
for all $i \in \{1, \dots, n\}.$

Next, we present two error bounds; one for the GSH and one for the GCSH. These error bounds follow immediately from Theorems \ref{thm:mainprop}$(iii)$ and \ref{thm:ebcentered2}$(iii)$ respectively. 

\begin{corollary}[General error bounds for one row of a Hessian]
Let $f:\dom f \subseteq \R^n\to\R$ be $\mathcal{C}^{4}$ on $B_n(x^0;\overline{\Delta})$ where $x^0 \in \dom f$ is the point of interest and $\overline{\Delta}>0$. Denote by  $L_{\nabla^2 f}\geq 0$ and $L_{\nabla^3 f} \geq 0$ the Lipschitz constant of $\nabla^2 f$  and $\nabla^3 f$ on $\overline{B}_n(x^0;\overline{\Delta})$ respectively.  Let $S=he^i \in \R^n$ where $h \neq 0,$  and $\T  \in \R^{n \times k}.$  Assume that $B_n(x^0+he^i;\Delta_{T})\subset B_n(x^0;\overline{\Delta}).$  Then

\begin{enumerate}[(i)]
\item
\begin{align*}
\left \Vert  \Projuv \nabla_s^2 f(x^0;S,\T)- \Projuv \nabla^2 f(x^0)  \right \Vert&= \left \Vert  \nabla_s^2 f(x^0;S,\T)- \Projuv \diag(e^j)\nabla^2 f(x^0) \right \Vert  \\
&\leq    4\sqrt{k} L_{\nabla^2 f} \left ( \frac{\Delta_u}{\Delta_l} \right ) \Vert \widehat{\T}^\dagger \Vert \Delta_u,
 \end{align*}
 and
\item 
\begin{align*}
\left \Vert  \Projuv \cshesst - \Projuv \nabla^2 f(x^0)  \right \Vert&= \left \Vert  \cshesst - \Projuv \Diag(e^i)\nabla^2 f(x^0)  \right \Vert  \\
&\leq  2  \sqrt{k}L_{\nabla^3 f} \left ( \frac{\Delta_u}{\Delta_l} \right )   \Vert  \widehat{\T}^\dagger  \Vert \Delta^2_u.
\end{align*}
\end{enumerate}
\end{corollary}

 Note that $\Vert (\widehat{S}^\top)^\dagger \Vert$ does not appear in the previous error bounds since
$\Vert (\widehat{S}^\top)^\dagger \Vert=1.$
One simple choice to approximate all entries of  the $i$\textsuperscript{th} row  is to choose
\begin{align} 
S&=he^i,~\T=h {\Id}_n \label{eq:approxRowi}
\end{align}
where $h \neq 0.$ In this case,  $\nabla_s^2 f(x^0;S,\T)$ is an order-1 accurate approximation of  the whole $i$\textsuperscript{th} row of the  Hessian. This choice uses the set of sample points     $x^0,x^0+he^i,x^0\oplus h{\Id}_n$ and $x^0+he^i \oplus h{\Id}_n.$ 
In this case, the number of function evaluations is 
$$1+1+n+n-1=2n+1.$$
We subtract one in the previous equation since one point is  reused: $x^0+he^i.$
Note that $2n+1 \leq (n+1)(n+2)/2$  for all $n \in \{1, 2, \dots\}.$ Therefore, if we are only interested by the entries of row/column $i$, setting $S$ and $\T$ in this fashion saves function evaluations  compared to using a minimal poised set for GSH. 

\begin{figure}[h]
\centering
\includegraphics[scale=0.4]{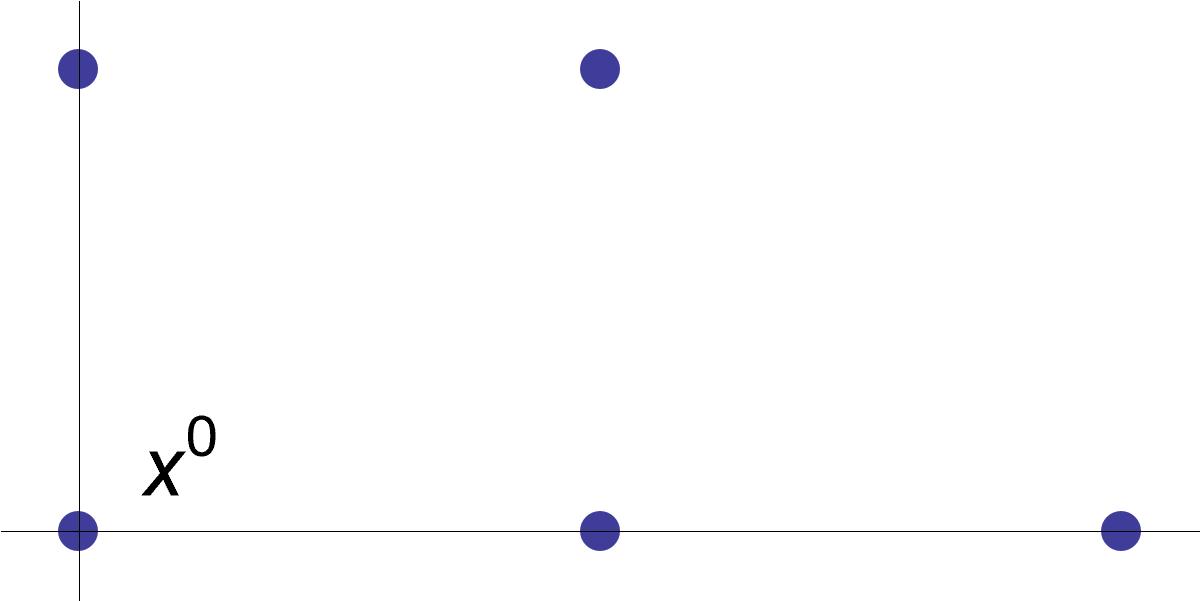}
\caption{Sample set  created while computing the GSH with  $S=e^1_2,\T=\Id_2$ as described in Eq. \eqref{eq:approxRowi}}
\end{figure} 

To obtain an order-2 accurate approximation of the whole $i$\textsuperscript{th} row/column of the Hessian, we may choose once again $S$ and $\T$ as defined in Equation \eqref{eq:approxRowi}. In this case, the  set of sample points is $x^0,x^0 \pm he^i,x^0\oplus (\pm h{\Id}_n), x^0+he^i \oplus h{\Id},$ and $x^0-he^i \oplus -{\Id}_n.$   Two sample points are reused: $x^0 \pm he^i.$ The number of function evaluations is 
$$1+2(2n)=4n+1.$$
Note that $4n+1 < n^2+n+1$ when $n \geq 4.$ Therefore, if  $n \in \{1, 2, 3\},$ then  using  a  minimal poised set for GCSH is preferable since it uses fewer function evaluations and provides an approximation of the full Hessian.

\begin{figure}[h]
\centering
\includegraphics[scale=0.4]{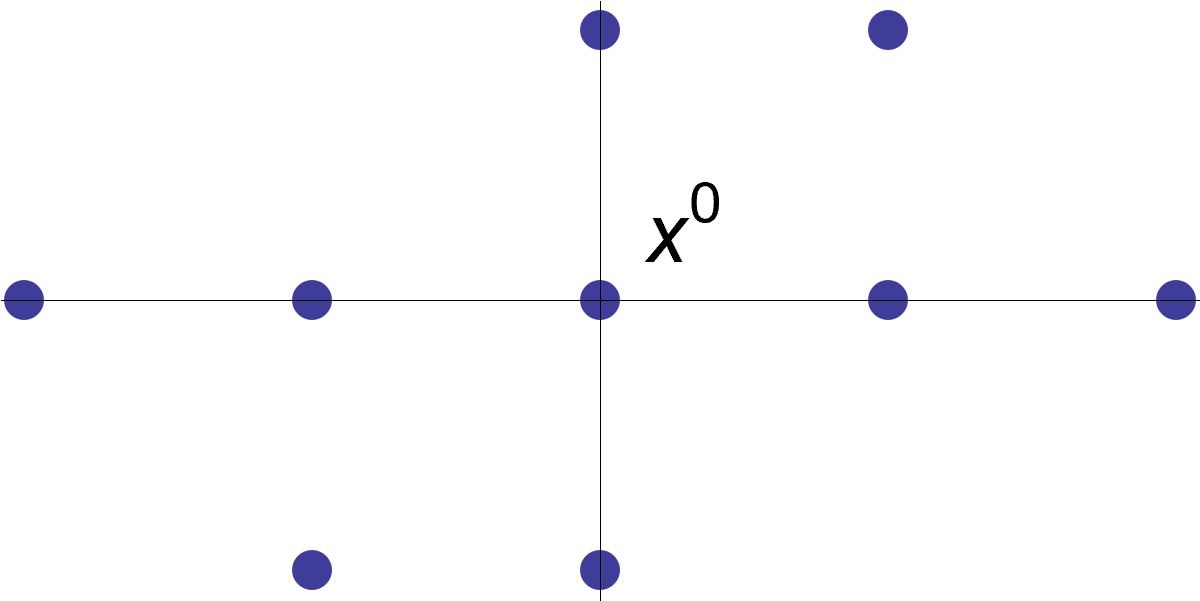}
\caption{Sample set  created while computing the GCSH with  $S=e^1_2,\T=\Id_2$ as described in Eq. \eqref{eq:approxRowi}}
\end{figure}

It seems that the minimum number of function evaluations to obtain an order-1 accurate approximation of a full row/column in a Hessian is $2n+1.$ To obtain an order-2 accurate approximation of a full row/column, the minimum number seems to be $n^2+n+1$ when $n \in \{1, 2, 3\}$ and $4n+1$ when $n \geq 4.$ Future research could focus on mathematically proving or disproving this claim.

In the following section, we investigate how to approximate a Hessian-vector product with the GSH and the GCSH.

\section{Approximating a Hessian-vector product}\label{sec:hvp}
In this section, we discuss how to approximate a Hessian-vector product (HVP) using the GSH and the GCSH.  Error bounds are provided for both techniques. The number of function evaluations required to accurately approximate a HVP is provided. We show that the GSH can provide an order-1 accurate approximation of a HVP  with $2n+1$ function evaluations and  the GCSH can provide an order-2 accurate approximation  with  $4n-1$ function evaluations. We begin  by defining the projection of a vector onto a  matrix $S$ and then provide an error  bound for the GSH. 

Given a matrix $S \in \R^{n \times m},$ the projection of a vector $w \in \R^n$ onto $S$ is denoted by $\Proj_S w$ and defined by $$\Proj_S w=(S^\top)^\dagger S^\top w.$$
This projection operator makes it possible to define an error bound for the GSH/GCSH in the case where $S$ is not full row rank (this projection operator has been previously defined in \cite{hare2020error}). In this case, it is only possible to define an error bound between the GSH/GCSH-vector product and some of the entries of  the HVP $\nabla^2 f(x^0) v.$ The appropriate entries are obtained by considering $\Proj_S (\nabla^2 f(x^0)v)$.

\begin{theorem}[Error bound for the GSH  when approximating a HVP]\label{thm:hvpeb}
Let $f:\dom f \subseteq \R^n\to\R$ be $\mathcal{C}^{3}$ on $B_n(x^0;\overline{\Delta})$ where $x^0 \in \dom f$ is the point of interest and $\overline{\Delta}>0$. Denote by $L_{\nabla^2 f}\geq 0$ the Lipschitz constant of $\nabla^2 f$ on $\overline{B}_n(x^0;\overline{\Delta})$. Let $v$ be a non-zero vector in $\R^n.$ Let $S=\begin{bmatrix} s^1&s^2&\cdots&s^m\end{bmatrix} \in \R^{n \times m}$ and $\T=hv,$ where $h \neq 0$. Assume that $B_n(x^0;\Delta_{\T})\subset B_n(x^0;\overline{\Delta})$ and $B_n(x^0+s^j;\Delta_{\T})\subset B_n(x^0;\overline{\Delta})$  for all $j \in \{1, \dots m\}.$ Then 
\begin{align*}
    \Vert \Proj_S (\nabla_s^2 f(x^0;S,\T)v)-\Proj_S(\nabla^2 f(x^0)v)\Vert &= \Vert \nabla_s^2 f(x^0;S,\T)v)-\Proj_S(\nabla^2 f(x^0)v)\Vert \leq  4\sqrt{m} L_{\nabla^2 f} \frac{\Delta_u}{\Delta_\ell} \Vert (\widehat{S}^\top)^\dagger  \Vert \Vert v \Vert  \Delta_u. 
\end{align*}
\end{theorem}
\begin{proof}
    From Definition \ref{def:mpinv}(ii), we know $(S^\top)^\dagger S^\top (S^\top)^\dagger=(S^\top)^\dagger.$  Hence, we obtain the first equality. To make notation more compact, let $H=\nabla^2 f(x^0).$ 
 We have 
\begin{align}
\Vert \shessv v- \Proj_S (Hv) \Vert &= \Vert (S^\top)^\dagger \delta_s^2 f(x^0;S,hv) v-(S^\top)^\dagger S^\top Hv \Vert \notag \\
&\leq \frac{\Vert (\widehat{S}^\top)^\dagger \Vert}{\Delta_S}  \Vert \delta_s^2 f(x^0;S,hv) v-S^\top Hv\Vert \notag \\
&=\frac{\Vert (\widehat{S}^\top)^\dagger \Vert}{\Delta_S}  \Vert \sqrt{\sum_{j=1}^m \vert (\nabla_s f(x^0+s^j;hv)-\nabla_s f(x^0;hv))^\top v-(s^j)^\top Hv\vert^2} \notag \\
&=\medmath{\frac{\Vert (\widehat{S}^\top)^\dagger \Vert}{\Delta_S}  \sqrt{\sum_{j=1}^m \vert \left ( (hv^\top)^\dagger \delta_s f(x^0+s^j;hv)-(hv^\top)^\dagger \delta_s f(x^0;hv)\right )^\top v-(s^j)^\top Hv\vert^2}} \notag \\
&=\frac{\Vert (\widehat{S}^\top)^\dagger \Vert}{\Delta_S}  \sqrt{\sum_{j=1}^m \left \vert  \frac{1}{h}\left (\delta_s f(x^0+s^j;hv)-\delta_s f(x^0;hv) \right )-(s^j)^\top Hv \right \vert^2} \notag \\
&=\medmath{\frac{\Vert (\widehat{S}^\top)^\dagger \Vert}{\vert h \vert \Delta_S}  \sqrt{\sum_{j=1}^m \vert f(x^0+s^j+hv)-f(x^0+s^j)-f(x^0+hv)+f(x^0)-(s^j)^\top H(hv)\vert^2}.}  \label{eq:plugin}
\end{align}
Each of the function values $f(x^0+s^j+hv), f(x^0+s^j)$ and $f(x^0+hv)$ may be written as a second-order Taylor expansion about $x^0$ plus a remainder term $R_2(x^0; \cdot).$ Hence, for all $j \in \{1, \dots, m\},$ we have
\begin{align}
&\left \vert f(x^0+s^j+hv)-f(x^0+s^j)-f(x^0+hv)+f(x^0)-(s^j)^\top H(hv) \right \vert \notag \\ 
& \quad \leq \vert R_2(x^0;s^j+hv) \vert +\vert R_2(x^0;s^j) \vert +\vert R_2(x^0;hv) \vert \notag \\
& \quad \leq \ \frac{1}{6}L_{\nabla^2 f} \Vert s^j+hv \Vert^3+\frac{1}{6}L_{\nabla^2 f} \Vert s^j\Vert^3+\frac{1}{6}L_{\nabla^2 f} \Vert hv \Vert^3  \notag \\
&\leq \frac{1}{2}L_{\nabla^2 f}(\Delta_S+\Delta_{\T} )^3. \label{eq:bound1} 
\end{align}
Substituting \eqref{eq:bound1} in \eqref{eq:plugin}, we obtain
\begin{align*}
    \Vert \shessv v- \Proj_S (Hv) \Vert &\leq  \frac{1}{\vert h \vert \Delta_S} \Vert (\widehat{S}^\top)^\dagger \Vert \sqrt{\sum_{j=1}^m \left ( \frac{1}{2} L_{\nabla^2 f}(\Delta_S+\Delta_{\T} )^3 \right )^2} \\
    &\leq 4\sqrt{m} L_{\nabla^2 f} \frac{\Delta_u}{\Delta_\ell}  \Vert (\widehat{S}^\top)^\dagger  \Vert \Vert v \Vert \Delta_u. 
\end{align*} \qed
\end{proof}

The previous theorem shows that we obtain  an order-1 accurate approximation of  $\nabla^2 f(x^0)v$  by taking $S$ to be  full row rank and $\T=hv, h \neq 0.$ If $f$ is a polynomial of degree 2 or less,   then we obtain a perfectly accurate approximation of the HVP.  If $S$ is not full row rank, then we obtain an order-1 accurate approximation of the HVP for some of the entries in the vector $\nabla^2f(x^0)v.$ The term $\Delta_u/\Delta_\ell$ suggests that if the radii are decreased, then they should be decreased by the same ratio  to keep the term  $\Delta_u/\Delta_\ell$ constant.

Taking $S=\begin{bmatrix} s^1&\cdots&s^n \end{bmatrix} \in \R^{n \times n}$, the GSH $\nabla_s^2 f(x^0;S,hv)$ produces the sample points  $$x^0, x^0+hv, x^0 \oplus S, \quad \text{and} \quad x^0+hv \oplus S.$$ 
This gives $2n+2$ function evaluations in the case where all points are distinct. To save one function evaluation, we may take $s^j=\pm hv$ for one $j \in \{1, \dots, n\}$. Hence, by taking $S \in \R^{n \times n}$ to be full rank and such that  $s^j=\pm hv$ for one $j \in \{1, \dots, n\},$ we obtain an order-1 accurate approximation of the HVP $\nabla^2 f(x^0)v$  and it costs $2n+1$ function evaluations. This is  one order of magnitude less than approximating a full Hessian with the GSH using a minimal poised set which requires  $(n+1)(n+2)/2.$ Note that $2n+1<(n+1)(n+2)/2$ whenever $n \geq 2.$

\begin{figure}[ht]
\centering
\includegraphics[scale=0.4]{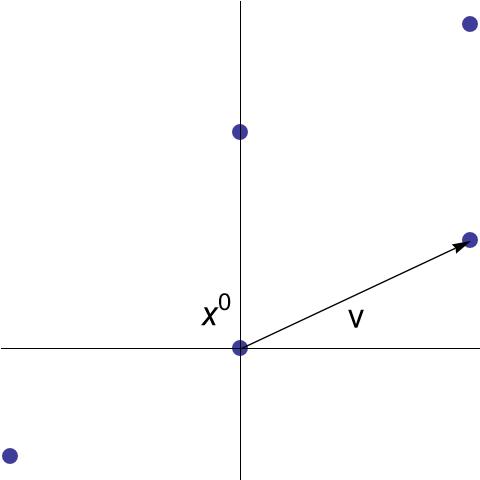}
\caption{Sample set  created while computing $\nabla_s^2 f(x^0;S,v)$ with $S=\begin{bmatrix} -v&e^2_2\end{bmatrix}$}
\end{figure}     

We now provide an error bound for the GCSH. The error bound shows that the GCSH is an order-2 accurate approximation of $\nabla^2 f(x^0)v$  under some assumptions.

\begin{theorem}[Error bound for the GCSH when approximating HVP]\label{thm:ebcenteredHVP}
Let $f:\dom f \subseteq \R^n\to\R$ be $\mathcal{C}^{4}$ on $B_n(x^0;\overline{\Delta})$ where $x^0 \in \dom f$ is the point of interest and $\overline{\Delta}>0$. Denote by $L_{\nabla^3 f}$ the Lipschitz constant of $\nabla^3 f$ on $\overline{B}_n(x^0;\overline{\Delta}).$ Let $v$ be a non-zero vector in $\R^n,$ $S=\begin{bmatrix} s^1&s^2&\cdots&s^m\end{bmatrix} \in \R^{n \times m},$ and $\T=hv$ where $h \neq 0.$ Assume that the ball $B_n(x^0 + s^j;\Delta_{\T})\subset B_n(x^0;\overline{\Delta})$ for all $j.$  Then  
\begin{align}
&\left \|\Proj_S (\cshesst v)- \Proj_S (\nabla^2 f(x^0)v) \right \|=\left \|  \cshesst- \Proj_S (\nabla^2 f(x^0)v) \right \| \leq 2  \sqrt{m}L_{\nabla^3 f} \frac{\Delta_u}{\Delta_l} \left \Vert  (\widehat{S}^\top)^\dagger  \right \Vert \Vert v \Vert  \Delta^2_u. \label{eq:hvpgcsh} 
\end{align}
\end{theorem}

\begin{proof}
The equality follows from Definition \ref{def:mpinv}(ii). To make notation more compact, let $H=\nabla^2 f(x^0)$. We have 
\small\begin{align}
&\left \|  \cshessv- \Proj_S (Hv) \right \| \notag \\
&=\left \Vert (S^\top)^\dagger \delta_c^2(x^0;S,hv)v-(S^\top)^\dagger S^\top Hv \right \Vert \notag \\
&=\frac{\Vert (\widehat{S}^\top)^\dagger \Vert}{\Delta_S} \left \Vert \delta_c^2(x^0;S,hv)v- S^\top Hv \right \Vert \notag \\
&=\frac{\Vert (\widehat{S}^\top)^\dagger \Vert}{\Delta_S} \sqrt{\sum_{j=1}^m \left \vert \frac{1}{2}\left ( \delta_s f(x^0+s^j;hv)+\delta_s f(x^0-s^j;-hv)-\delta_s f(x^0;hv)-\delta_s f(x^0;-hv) \right )(hv)^\dagger v-(s^j)^\top Hv \right \vert^2} \notag \\
&=\medmath{\frac{\Vert (\widehat{S}^\top)^\dagger \Vert}{\vert h \vert \Delta_S} \sqrt{\sum_{j=1}^m \left \vert \frac{1}{2} \left (f(x^0+s^j+hv)-f(x^0+s^j)+f(x^0-s^j-hv)-f(x^0-s^j)-f(x^0+hv)-f(x^0-hv)+2f(x^0) \right )-(s^j)^\top H(hv) \right \vert^2}.}  \label{eq:plugin2}
\end{align}\normalsize
Each of the function values $f(x^0+s^j+hv), f(x^0+s^j), f(x^0-s^j-hv), f(x^0-s^j), f(x^0+hv)$ and $f(x^0-hv)$ may be written as a third-order Taylor expansion about $x^0$ plus a remainder term $R_3(x^0;\cdot).$ It follows that 
\begin{align}
&\medmath{\left \vert \frac{1}{2} \left (f(x^0+s^j+hv)-f(x^0+s^j)+f(x^0-s^j+hv)-f(x^0-s^j)-f(x^0+hv)-f(x^0-hv)+2f(x^0) \right )-(s^j)^\top H(hv) \right \vert} \notag\\
&\leq \frac{1}{2} \left (\vert R_3(x^0;s^j+hv)\vert +\vert R_3(x^0;s^j)\vert + \vert R_3(x^0;+s^j+hv)\vert +\vert R_3(x^0;-s^j) \vert +\vert R_3(x^0;hv \vert +\vert R_3(x^0;-hv)\vert \right ) \notag \\
     &\leq \frac1{48}  \left (L_{\nabla^3 f} \Vert s^j+hv\Vert^4+ L_{\nabla^3 f} \Vert s^j \Vert^4+L_{\nabla^3 f} \Vert -s^j+hv\Vert^4+L_{\nabla^3 f} \Vert -s^j \Vert^4 +L_{\nabla^3 f} \Vert hv \Vert^4 +L_{\nabla^3 f} \Vert -hv \Vert^4   \right )\notag\\
     &\leq \frac{1}{8} L_{\nabla^3 f} (\Delta_{S}+\Delta_{\T})^4. \label{eq:remainder}
\end{align}
Substituting the bound \eqref{eq:remainder} in \eqref{eq:plugin2}, we get 
\begin{align*}
\left \|  \cshesst- \Proj_S (Hv) \right \| &\leq   \frac{\sqrt{m}}{8}  L_{\nabla^3 f} \frac{(\Delta_S+\Delta_{\T})^4}{\Delta_S \Delta_{\T}} \left \Vert(\widehat{S}^\top)^\dagger \right  \Vert \Vert v \Vert  \\
&\leq 2 \sqrt{m}  L_{\nabla^3 f} \frac{\Delta_u}{\Delta_l}  \left \Vert(\widehat{S}^\top)^\dagger \right \Vert \Vert v \Vert \Delta_u^2.
\end{align*}
\qed
\end{proof}

The previous error bound shows that the GCSH provides an order-2 accurate approximation of  $\nabla^2 f(x^0)v$  whenever $S$ is full row rank and $\T=hv, h \neq 0.$ If $f$ is a polynomial of degree 3 or less, then it provides a perfectly accurate approximation of $\nabla^2f(x^0)v$.

Taking $S=\begin{bmatrix} s^1&\cdots&s^n\end{bmatrix} \in \R^{n \times n}$  and $\T=hv,$  the GCSH $\nabla^2_c f(x^0;S,hv)$ produces the sample points $$x^0, x^0 \pm hv, x^0 \oplus (\pm S), x^0+hv \oplus S, x^0-hv \oplus (-S).$$ 
Assuming all points are distinct, $4n+3$ function evaluations are required.  By taking $s^j=-hv, h>0,$ for one $j \in \{1, \dots, n\},$  4 function evaluations are saved, and we get  an order-2 accurate approximation of the HVP for $4n-1$ function evaluations. This number is one order of magnitude less than the number of function evaluations required when approximating  a full Hessian with the GCSH using a minimal poised set ($n^2+n+1$).  Note that $4n-1 < n^2+n+1$ whenever $n \geq 3.$

\begin{figure}[ht]
\centering
\includegraphics[scale=0.4]{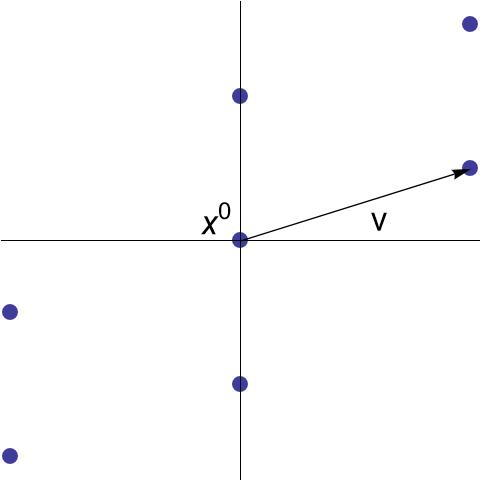}
\caption{Sample set  created while computing $\nabla_c^2 f(x^0;S,v)$ with $S=\begin{bmatrix} -v&e^2_2\end{bmatrix}$}
\end{figure}     

In the next section, we investigate how to approximate higher-order derivatives by defining a technique  with a similar structure than the  generalized simplex gradient and the GSH.
                         
\FloatBarrier

\section{Approximating order-P derivatives} \label{sec:nablaP}
Now that we have a general method to approximate first-order derivatives called the generalized simplex gradient and a general method to approximate second-order derivatives called the generalized simplex Hessian, we may develop a general method to approximate  $P$-order derivatives.  The object containing  all $P$-order derivatives  can  be viewed as a $P$-dimensional matrix. We begin  by providing a formula to approximate  all third-order derivatives  and then we propose a formula to compute $P$-order derivatives.

We refer to  $\nabla^3 f(x^0)$ as the \emph{Tressian} of $f$ at $x^0.$ A Tressian can be viewed as a three-dimensional matrix (a tensor) $\mathbf{M}$ in $\R^{n \times n \times n}$ where the third dimension represents the depth of  $\mathbf{M}$.  In this section, a tensor $\mathbf{M} \in \R^{r \times c \times p}$ will be either thought of as an object  containing $r$ floors where each floor is a matrix in $\R^{c \times p}$ or as an object containing $p$ layers where each layer is a matrix in $\R^{r\times c}.$ 
A  tensor $\mathbf{M} \in \R^{r \times c \times p}$ is written by floor in the following way:

\begin{align} \label{eq:3dmatrixbyfloor}
    \mathbf{M}&=\begin{bmatrix} F_1\\F_2 \\ \vdots \\F_r \end{bmatrix}_{(r,\cdot,\cdot)} 
\end{align}
where $F_i \in \R^{c \times p}$ for all $i \in \{1, 2, \dots, r\}$ and $[F_i]_{j,k}=\mathbf{M}_{i,j,k}$ for all $i,j,k.$  The subscript in \eqref{eq:3dmatrixbyfloor} is used to make it clear that $\mathbf{M}$ is written by floor.  The tensor $\mathbf{M}$ can also be written in terms of layers:
\begin{align*}\label{eq:byFloor}
    \mathbf{M}&=\begin{bmatrix} L_1\\L_2\\ \vdots \\ L_p \end{bmatrix}_{(\cdot,\cdot,p)} 
\end{align*}
where $[L_k]_{i,j}=\mathbf{M}_{i,j,k}$ for all $i,j,k.$

We are now ready to introduce the formula to approximate  $\nabla^3 f(x^0).$ The technique requires one more set  of matrices of directions than the generalized simplex Hessian. The letter $U$ is used to denote this new  of set of matrices. This set of matrices associated could contain  $k_1+k_2+\dots+k_m$ different matrices.  To keep things relatively simple, we provide the formula for the case where all matrices $T_j,$ are equal and all matrices $U_k \in \R^{n \times \ell_k}$ are equal. To emphasize this special case where all matrices $U_k$ are equal, we use the notation $\overline{U} \in \R^{n \times \ell}.$ Hence, the three matrices of directions involved in the computation of the approximation technique are $S \in \R^{n \times m}, \T \in \R^{n \times k},$ and $\overline{U} \in \R^{n \times \ell}.$  

Before introducing the approximation technique, we define the multiplication of a two-dimensional matrix with a tensor.

Let $A \in \R^{n \times m}$ and let $\mathbf{M} \in \R^{m \times n \times p}.$ Let $\mathbf{M}$ be written as layers: 
$$\mathbf{M} = \begin{bmatrix} L_1\\L_2\\\vdots\\L_p \end{bmatrix}_{(\cdot,\cdot,p)} \in \R^{m \times n \times p},$$
where $[L_k]_{i,j}=\mathbf{M}_{i,j,k}.$ 
Then 
$$A \otimes \mathbf{M}= \begin{bmatrix} AL_1\\AL_2\\\vdots\\AL_p \end{bmatrix}_{(\cdot,\cdot,p)} \in \R^{n \times n \times p}$$

\begin{definition}[Generalized simplex Tressian]
    Let $f:\dom f \subseteq \R^n \to \R$ and let $x^0 \in \dom f$ be the point of interest.  Let $S=\begin{bmatrix} s^1&s^2&\cdots &s^m \end{bmatrix}  \in \R^{n \times m}$ and $ \T \in \R^{n \times k}, \overline{U} \in \R^{n \times \ell}$ with  the set of sample points $\mathcal{S}(x^0;S,\overline{T},\overline{U})$ contained in $ \dom f.$  The \emph{generalized simplex Tressian}  of $f$ at $x^0$ over $S, \overline{T}$ and $\U$ is denoted by $\stress$ and defined by
\begin{equation*}\stress =(S^\top)^\dagger \otimes \dst f(x^0;S,\T,\U) \in \R^{n \times n \times n},\end{equation*}
where
\begin{equation*}
    \dst f(x^0;S,\T,\U)=\left[\begin{array}{c}(\ns^2 f(x^0+s^1;\T,\U)-\ns^2 f(x^0;\T,\U))^\top\\(\ns^2 f(x^0+s^2;\T,\U)-\ns^2 f(x^0;\T,\U))^\top\\\vdots\\(\ns^2 f(x^0+s^m;\T,\U)-\ns f(x^0;\T,\U))^\top\end{array}\right]_{(m,\cdot,\cdot)} \in \R^{m \times n \times n}.
\end{equation*}
\end{definition}

Recursively,  we may now define a simple formula to approximate order-$P$ derivatives of a function at a point of interest $x^0 \in \R^n$.  Before introducing the formula,  notation  needs to be slightly modified  to make it easier to discuss   general order-$P$ derivatives.
To approximate  order-$P$ derivatives, we use a matrix $S_1 \in \R^{n \times m_1},$ and set of matrices $S_2, S_3, \dots, S_P.$ To keep notation relatively simple, we consider the case  where all matrices of directions are the same in the sets $S_2, \dots, S_P.$ As before, we write  $\overline{S_i}$ to emphasize that all matrices of directions are  identical in each set $\overline{S_i}, i \in \{2, 3, \dots, P\}.$ A matrix in the set $\overline{S_i}$ has dimensions $n \times m_i,$ for $i \in \{2, 3, \dots, P\}.$

The transpose of a $P$-dimensional matrix $\mathbf{M} \in \R^{n\times m_1 \times \dots \times m_{P-1}}$  is denoted by $\mathbf{M}^\top$ where the entries of $\mathbf{M}^\top$ are equal to
$$[\mathbf{M}^\top]_{i,j_1, \dots,j_{P-1}}=\mathbf{M}_{j_{P-1}, j_{P-2},\dots, j_1,i}, \quad i \in \{1, 2, \dots, n\}, j_k \in \{1, 2, \dots, m_k\}, k \in \{1, 2, \dots, P-1\}.$$

\begin{definition}[Order-$P$ simplex derivative matrix]
    Let $f:\dom f \subseteq \R^n \to \R$ and let $x^0 \in \dom f$ be the point of interest.  Let $S_1=\begin{bmatrix} s^1&s^2&\cdots &s^{m_1} \end{bmatrix}  \in \R^{n \times m_1}$ and $ \overline{S_i} \in \R^{n \times m_i}$ for all $i \in \{2, 3, \dots, P\}$  with  all sample points contained in $ \dom f.$  The \emph{order-$P$ simplex derivative tensor of $f$ at $x^0$ over $S_1, \overline{S_2}, \dots, \overline{S_P}$}  is denoted by $\nabla_s^P f(x^0;S_1,\overline{S_2}, \dots, \overline{S_P})$ and defined by
\begin{equation*}
\nabla_s^P f(x^0;S_1,\overline{S_2},\dots, \overline{S_P}) =(S_1^\top)^\dagger \otimes \delta_s^P f(x^0;S_1,\overline{S_2},\dots, \overline{S_P}) \in \R^{n \times n \times \dots \times n},\end{equation*}
where
\begin{equation*}
    \delta_s^Pf(x^0;S_1,\overline{S_2},\dots,\overline{S_P})=\left[\begin{array}{c}(\ns^{P-1} f(x^0+s^1;\overline{S_2},\dots,\overline{S_P})-\ns^{P-1} f(x^0;\overline{S_2},\dots,\overline{S_P}))^\top\\(\ns^{P-1} f(x^0+s^2;\overline{S_2},\dots,\overline{S_P})-\ns^{P-1} f(x^0;\overline{S_2},\dots,\overline{S_P}))^\top\\\vdots\\(\ns^{P-1} f(x^0+s^{m_1};\overline{S_2},\dots,\overline{S_P})-\ns^{P-1} f(x^0;\overline{S_2},\dots,\overline{S_P}))^\top\end{array}\right]_{(m_1,\cdot,\dots,\cdot)}.
\end{equation*}
\end{definition}

The previous definition provides a compact formula that can be  (relatively) easily implemented in a software such as MATLAB.
\section{Conclusion} \label{sec:conclusion}


In Section \ref{sec:partialHessian}, we provided details on how  to choose the matrices $S$ and $T_j$ when we are only interested in a proper subset of the entries of the Hessian.  In particular, we investigated how to approximate the diagonal entries of a  Hessian, the off-diagonal entries of a Hessian, and a row/column  of a Hessian. The number of function evaluations  to obtain an order-1 accurate approximation, or an order-2 accurate approximation of the entries of the Hessian of interest has been discussed. This shows that the GSH is a valuable tool  to approximate either a full Hessian or a proper subset of the entries of a Hessian. In both cases,  the error bounds provided show that the error can be controlled by the individual conducting the optimization process. The GSH and the GCSH are simple approximation techniques based on matrix algebra that can be used to do ``everything'' related to Hessian approximation. These techniques can be  easily implemented in a software such as MATLAB.

The relation between the CSHD introduced in \cite{jarry2022approximating} and the GCSH is clarified. It is shown that the CSHD is equal to the GCSH whenever $S$ is a partial diagonal matrix with full column rank and $T_j=-s^j$ for all $j$ (Theorem \ref{thm:shequalcshd}).

In Section \ref{sec:hvp}, it has been shown how  the GSH and the GCSH can be used to approximate a Hessian-vector product $\nabla^2 f(x^0)v.$ Essentially, the ``trick'' is to choose the second matrix of directions $\T$ to be equal to a non-zero multiple of the vector $v.$ 

In Section \ref{sec:nablaP}, the approximation  technique is generalized to  higher-order derivatives. First,  it is discussed how to obtain  an approximation of the third-order derivatives.   Then a simple recursive formula is introduced to computer order-$P$ derivatives of a function at a point of interest. 

It remains to verify if an order-1 accurate approximation of the main diagonal of a Tressian can be obtained for free in terms of function evaluations, if an order-2 accurate approximation of the Hessian has been previously computed via the GCSH. It is reasonable to believe that it is the case since an order-2 accurate approximation of the gradient using the generalized centered simplex gradient provides an order-1 accurate approximation  of the diagonal of the Hessian for free.

On a final note, an implementation in MATLAB of each approximation technique discussed in this paper is available upon request.
\newpage 

\vspace{0.5cm}
\noindent \small{\textbf{Author contributions} All authors contributed equally to this work.}

\vspace{0.5cm}

\noindent \small{\textbf{Funding} Chayne Planiden was supported by the University of Wollongong.}

\vspace{0.5cm}

\noindent \small{\textbf{Data availability} No datasets were generated or analyzed during the current study.}

\vspace{0.5cm}

\noindent \small{\textbf{Code availability} The codes are available upon request.}

\section*{Declarations}

\noindent \small{\textbf{Ethical approval} The authors declare that they followed all the rules of a good scientific practice.}

\vspace{0.5cm}

\noindent \small{\textbf{Consent to participate} All authors approve their participation in this work.}

\vspace{0.5cm}

\noindent \small{\textbf{Consent for publication} The authors approve the publication of this research.}

\vspace{0.5cm}

\noindent \small{\textbf{Conflict of interest} The authors declare no competing interests.}

\vspace{0.5cm}

\noindent \small{\textbf{Human and animal ethics} Not applicable.}
\clearpage

\normalsize
\bibliographystyle{siam}
\bibliography{bibliography}
\end{document}